
\def\LaTeX{%
  \let\Begin\begin
  \let\End\end
  \let\salta\relax
  \let\finqui\relax
  \let\futuro\relax}

\def\UK{\def\our{our}\let\sz s}
\def\USA{\def\our{or}\let\sz z}

\UK



\LaTeX

\USA


\salta

\documentclass[twoside,12pt]{article}
\setlength{\textheight}{24cm}
\setlength{\textwidth}{16cm}
\setlength{\oddsidemargin}{2mm}
\setlength{\evensidemargin}{2mm}
\setlength{\topmargin}{-15mm}
\parskip2mm


\usepackage[usenames,dvipsnames]{color}
\usepackage{amsmath}
\usepackage{amsthm}
\usepackage{amssymb}
\usepackage[mathcal]{euscript}
\usepackage{enumitem}
%
%
\usepackage{cite}
%
%
%


\definecolor{viola}{rgb}{0.3,0,0.7}
\definecolor{ciclamino}{rgb}{0.5,0,0.5}
\definecolor{rosso}{rgb}{0.8,0,0}

\def\andrea #1{{\color{red}#1}}
\def\dafare #1{{\color{red}#1}}
\def\cit #1{{\color{OliveGreen}#1}}

\def\andrea #1{#1}
\def\dafare  #1{#1}
\def\cit  #1{#1}


\bibliographystyle{plain}


%

\finqui

\def\Beq{\Begin{equation}}
\def\Eeq{\End{equation}}
\def\Bsist{\Begin{eqnarray}}
\def\Esist{\End{eqnarray}}

\def\Bthm{\Begin{theorem}}
\def\Ethm{\End{theorem}}

\def\Bcor{\Begin{corollary}}
\def\Ecor{\End{corollary}}
\def\Brem{\Begin{remark}\rm}
\def\Erem{\End{remark}}

\def\Bcenter{\Begin{center}}
\def\Ecenter{\End{center}}
\let\non\nonumber




\def\step #1 \par{\medskip\noindent{\bf #1.}\quad}


\def\Lip{Lip\-schitz}
\def\Holder{H\"older}
\def\Frechet{Fr\'echet}

\def\aand{\quad\hbox{and}\quad}

\def\wk{well-known}
\def\socal{so-called}
\def\lhs{left-hand side}
\def\rhs{right-hand side}
\def\sfw{straightforward}

\def\CH{Cahn--Hilliard}


\def\generaliz{generali\sz}


\def\multibold #1{\def\arg{#1}%
  \ifx\arg\pto \let\next\relax
  \else
  \def\next{\expandafter
    \def\csname #1#1#1\endcsname{{\bf #1}}%
    \multibold}%
  \fi \next}

\def\pto{.}

\def\multical #1{\def\arg{#1}%
  \ifx\arg\pto \let\next\relax
  \else
  \def\next{\expandafter
    \def\csname cal#1\endcsname{{\cal #1}}%
    \multical}%
  \fi \next}


\def\multimathop #1 {\def\arg{#1}%
  \ifx\arg\pto \let\next\relax
  \else
  \def\next{\expandafter
    \def\csname #1\endcsname{\mathop{\rm #1}\nolimits}%
    \multimathop}%
  \fi \next}

\multibold
qwertyuiopasdfghjklzxcvbnmQWERTYUIOPASDFGHJKLZXCVBNM.

\multical
QWERTYUIOPASDFGHJKLZXCVBNM.

\multimathop
ad dist div dom meas sign supp .


\def\accorpa #1#2{\eqref{#1}--\eqref{#2}}
\def\Accorpa #1#2 #3 {\gdef #1{\eqref{#2}--\eqref{#3}}%
  \wlog{}\wlog{\string #1 -> #2 - #3}\wlog{}}


\def\tonde #1{\left(#1\right)}

\def\graffe #1{\mathopen\{#1\mathclose\}}

\def\<#1>{\mathopen\langle #1\mathclose\rangle}
\def\norma #1{\mathopen \| #1\mathclose \|}

\def\iot {\int_0^t}

\def\intQt{\int_{Q_t}}
\def\intQtT{\int_{Q_t^T}}
\def\intQ{\int_Q}
\def\iO{\int_\Omega}

\def\dt{\partial_t}
\def\dn{\partial_n}
\def\dtt{\partial_{tt}}

\def\checkmmode #1{\relax\ifmmode\hbox{#1}\else{#1}\fi}

\def\aeQ{\checkmmode{a.e.\ in~$Q$}}

\def\aaQ{\checkmmode{for a.a.~$(x,t)\in Q$}}

\def\aat{\checkmmode{for a.a.~$t\in(0,T)$}}


\def\erre{{\mathbb{R}}}




\def\genspazio #1#2#3#4#5{#1^{#2}(#5,#4;#3)}
\def\spazio #1#2#3{\genspazio {#1}{#2}{#3}T0}

\def\L {\spazio L}
\def\H {\spazio H}
\def\W {\spazio W}

\def\C #1#2{C^{#1}([0,T];#2)}


\def\Lx #1{L^{#1}(\Omega)}
\def\Hx #1{H^{#1}(\Omega)}

\def\Luno{\Lx 1}
\def\Ldue{\Lx 2}
\def\Linfty{\Lx\infty}

\def\Huno{\Hx 1}
\def\Hdue{\Hx 2}



\let\theta\vartheta
\let\eps\varepsilon
\let\phi\varphi
\let\lam\lambda

\let\TeXchi\chi                         
\newbox\chibox
\setbox0 \hbox{\mathsurround0pt $\TeXchi$}
\setbox\chibox \hbox{\raise\dp0 \box 0 }
\def\chi{\copy\chibox}



\def\bQ{b_1}
\def\bO{b_2}
\def\bQh{b_3}
\def\bOh{b_4}
\def\bz{b_0}

\def\rmin{r_-}
\def\rmax{r_+}

\def\phQ{\phi_Q}

\def\phO{\phi_\Omega}

\def\mQ{\m_Q}
\def\mO{\m_\Omega}

\def\sO{\s_\Omega}

\def\sQ{\s_Q}

\def\Uad{\calU_{\ad}}

\def\vO{v^\Omega}

\def\Vp{V^*}

\def\normaV #1{\norma{#1}_V}
\def\normaH #1{\norma{#1}_H}

\def\normaVp #1{\norma{#1}_*}

\let\hat\widehat

\def\Pi{\hat\pi}


\def\cd{c_\delta}
\def\s{\sigma}  
\def\m{\mu}	    
\def\ph{\phi}	
\def\a{\alpha}	
\def\b{\beta}	
\def\d{\delta}  
\def\et{\eta}   
\def\th{\theta} 
\def\r{\rho}    
\def\bph{\bar\ph}  
\def\bm{\bar\m}    
\def\bs{\bar\s}    
\def\z{\zeta}      
\def\ch{\chi}      
\def\ps{\psi}      
\def\J{{\cal J}}
\def\S{{\cal S}}
\def\I2 #1{\int_{Q_t}|{#1}|^2}
\def\IT2 #1{\int_{Q_t^T}|{#1}|^2}
\def\IN2 #1{\int_{Q_t}|\nabla{#1}|^2}
\def\IO2 #1{\iO |{#1(t)}|^2}
\def\INO2 #1{\iO |\nabla{#1}(t)|^2}
\def\UR{{\cal U}_R}
\def\CP{{\bf (CP)}}
\def\A{{\cal A}}
\def\span{{\it span}}

\Begin{document}

\title{Optimal distributed control
		\\[0.3cm] 
		of an extended model of tumor growth
		\\[0.3cm] 
		with logarithmic potential}
\author{}
\date{}
\maketitle

\Bcenter
\vskip-1cm
{\large\sc Andrea Signori$^{(1)}$}\\
{\normalsize e-mail: {\tt andrea.signori02@universitadipavia.it}}\\[.25cm]
$^{(1)}$
\dafare{
{\small Dipartimento di Matematica e Applicazioni, Universit\`a di Milano--Bicocca}\\
{\small via Cozzi 55, 20125 Milano, Italy}
}
\Ecenter
\Begin{abstract}\noindent
This paper is intended to tackle the control problem associated with
an extended phase field system of \CH~type that is related to a tumor growth model.
This system has been investigated in previous contributions
from the viewpoint of well-posedness and asymptotic analyses.
Here, we aim to extend the mathematical studies around this system by introducing
a control variable and handling the corresponding control problem.
We try to keep the potential as general as possible, focusing our
investigation towards singular potentials, such as the logarithmic one.
We establish the existence of optimal control, the \Lip~continuity of
the control-to-state mapping and even its \Frechet~differentiability in suitable Banach spaces.
Moreover, we derive the first-order necessary conditions that an optimal control has to satisfy.

\vskip3mm
\noindent {\bf Key words:}
\dafare{Distributed optimal control, tumor growth, phase field model, Cahn--Hilliard equation, 
optimal control, necessary optimality conditions, adjoint system.}
\vskip3mm
\dafare{\noindent {\bf AMS (MOS) Subject Classification:} 35K61, 35Q92, 49J20, 49K20, 92C50.}
\End{abstract}

\vskip3mm

\pagestyle{myheadings}
\newcommand\testopari{\sc Signori}
\newcommand\testodispari{\sc Optimal distributed control of tumor growth model}
\markboth{\testodispari}{\testopari}

\salta
\finqui
\newpage
\section{Introduction}
\label{SEC_INTRO}
\setcounter{equation}{0}
In this paper, we deal with a distributed optimal control
problem for a system of partial differential equations
whose physical context is that of tumor growth dynamics. Our aim is to
devote this section to explain the general purpose of the work
and we postpone all the technicalities for the forthcoming
sections. In the next one, we will state precisely the problem 
and have the care to present in detail our notation and the mathematical framework
in which set the problem.
Here, let us only mention that with $\Omega \subset \erre^3$ we denote the 
set where the evolution takes place and, for a given final time $T>0$, we fix
\Beq
	\non
	Q:=\Omega \times (0,T) \quad \hbox{and} \quad \Sigma:=\Gamma \times (0,T).
\Eeq
The distributed control problem, referred as \CP, consists of minimizing
the \socal~cost functional
\Bsist
	&& \non
	\J (\ph, \s, u)   = 
	\frac \bQ 2 \norma{\ph - \phQ}_{L^2(Q)}^2
	+\frac \bO 2 \norma{\ph(T)-\phO}_{L^2(\Omega)}^2
	+ \frac \bQh 2 \norma{\s - \sQ}_{L^2(Q)}^2
	\\ && \quad
	+ \frac \bOh 2 \norma{\s(T)-\sO}_{L^2(\Omega)}^2
    + \frac \bz 2 \norma{u}_{L^2(Q)}^2,
    \label{costfunct}
\Esist
subject to the control contraints 
\Beq
    \label{Uad}
	u \in \Uad := \graffe{u \in L^{\infty}(Q): u_* \leq u \leq u^*\ \aeQ},
\Eeq
and to the state system
\Bsist
  & \a \dt \m + \dt \ph - \Delta \m = P(\phi) (\sigma - \mu)
  \quad \hbox{in $\, Q$}
  \label{EQprima}
  \\[0.2cm]
  & \mu = \beta \dt \phi - \Delta \phi + F'(\phi)
  \label{EQseconda}
  \quad \hbox{in $\,Q$}
  \\[0.2cm]
  & \dt \sigma - \Delta \sigma = - P(\phi) (\sigma - \mu) + u
  \label{EQterza}
  \quad \hbox{in $\,Q$}
  \\[0.4cm]
  & \dn \m = \dn \ph =\dn \s = 0
  \quad \hbox{on $\,\Sigma$}
  \label{BCEQ}
  \\[0.2cm]
  & \m(0)=\m_0,\, \ph(0)=\ph_0,\, \s(0)=\s_0
  \quad \hbox{in $\,\Omega.$}
  \label{ICEQ}
\Esist
\Accorpa\EQ EQprima ICEQ
Let us give just some overall indications on the involved quantities of the above equations.
The symbols $\bz, \bQ,\bO,\bQh,\bOh$ 
represent nonnegative constants, not all zero, while 
$\phQ, \phO, \sQ, \sO $,
$u_*$, and $u^*$ denote given functions.
As regards these latter, the first four model some targets, while 
the last two fix the box in which the control variable $u$ can be chosen.
Furthermore, $F$ and $P$ are nonlinearities, while \eqref{BCEQ} and
\eqref{ICEQ} are the boundary conditions and the initial conditions, 
respectively.

During the last decades, lots of models based on continuum mixture theory have been derived.
The above state system constitutes a variation on 
an approximation to a diffuse interface model for the dynamics 
of tumor growth proposed in \cite{HDZO} (see also \cite{HKNZ} and \cite{WZZ}), in which 
the velocity contributions are neglected and the attention is focused on the behavior of 
the state variables that model the fractions of the tumor cells and the 
nutrient-rich extracellular water, respectively. 
\dafare{Moreover, let us refer to \cit{\cite{GARL_1, GARL_2, GARL_3, GARL_4, GARLR, GAR}}, 
where transport mechanisms such as chemotaxis and active transport are also taken into account.}
Further investigations and mathematical models
related to biology can be found e.g. in \cite{DFRGM} and \cite{FLRS}.

Let us spend some words about the interpretation of the system \EQ,
and on the involved variables. 
The unknown $\ph$ is an order parameter which describes the tumor cell
fraction and assumes values between $-1$ and $+1$. These two extremes represent the pure phases,
say the tumor phase and the healthy cell phase, respectively. The second unknown $\m$ has the
interpretation, as usual for \CH~equation, of chemical potential and 
its relation with $\ph$ is precisely expressed by \eqref{EQseconda}.
The third unknown $\s$ consists of the nutrient-rich extracellular water volume fraction
and we assume that it takes values between $0$ and $1$ with the following property:
the closer to one, the richer of water the extracellular fraction is, while the closer to zero,
the poorer it is. As the nonlinearities are concerned,
we have that $F$ stands for a double-well potential, while 
$P$ models a proliferation function which we assume to be nonnegative and 
dependent on the phase variable.
To conclude the overview of the model, worth to
point out the different role of $\a$ and $\b$.
When $\a=0$, the equations \accorpa{EQprima}{EQseconda} becomes of viscous \CH~type or it
is pure \CH~equation depending on the fact that $\b$ is strictly positive or
vanishes, respectively. On the other hand, the presence of $\a$ gives to \eqref{EQprima}
a parabolic structure with respect to the variable $\m$.

As for the interpretation of the \CP~problem, our goal consists of 
finding a ``smarter" choice of 
$u\in\Uad$ such that, with its corresponding solution 
to \EQ, minimizes \eqref{costfunct}. 
Note that the control variable $u$ appears in \eqref{EQterza}, the equation describing
the nutrient evolution process. Thus, from the viewpoint of the model, it could 
represent a supply of a nutrient or a drug in chemotherapy. 
The cost functional we choose is a tracking-type one, namely we 
have fixed some a priori targets, say some a priori final configurations 
for the tumor cells and on the nutrient, 
and we try to find the control variable whose corresponding
solutions approximate better this fixed configuration.
Worth to insist on this fact: even if the better situation is the health of the patient,
our efforts are neither in the direction of minimizing the variable $\ph$, that has the meaning of
leading to the healthier configuration nor minimizing the variable $\s$ to reduce the tumor
expansion. In fact, we only 
try to handle the whole evolution process, acting on the choice of the control variable, 
to force a final configuration that for some practical reason should be desirable.
Obviously the ratios among the constants $\bz,\bQ,\bO,\bQh,\bOh$
implicitly describe which targets hold the leading part in our application.
To conclude the analysis, we focus our attention on the last term of \eqref{costfunct}.
From an abstract viewpoint, it represents the cost we have to pay to implement $u$, thus
in our framework it should be read as the rate of risks to afflict harm to the
patient by following that strategy. 
Finally, observe that we do not consider the cost functional to be
dependent on the chemical potential. Indeed,
from an interpretation point of view, we mainly care to handle the phase dynamics, and
it is not clear if including the variable $\m$ in the analysis is interesting for applications
(see the forthcoming Remark \ref{REM_modcost}).

At this general stage, let us perform a little overview of the literature.
The first systematic study on this system was carried out in \cite{CGH} and \cite{FGR},
where well-posedness and long-time behavior of the solutions were investigated
for a system very close to ours.
Moreover, quite recently, the system has 
been investigated with particular interest on the asymptotic analysis as the constants $\a$ 
and $\b$ go to zero. To this concern, we address to \cite{CGH}, \cite{CGRS_ASY}, and \cite{CGRS_VAN},
where the asymptotic analyses represent the core of the works. 
To our best knowledge, as the control theory is concerned, there are very few contributions
to this kind of system. In this regard, we refer to \cite{CGRS_OPT}, where
a control problem for a system without relaxation terms is performed. Even though 
we take inspiration from this work, the functional framework and the 
potentials setting significantly differ from ours. 
Nevertheless, the control theory related to different phase-field models based on 
the \CH~equation presents more contributions. Among others, we mention
\cite{CGMR_cons,CGMR_sing,CGS_nonst,CGS14,CS}.
Furthermore, since particular attention is devoted to singular potentials,
we point out \cite{CGS,GiMiSchi,MirZelik} and the vast 
list of references therein.

To conclude, let us sketch an outline of the work.
The first section is devoted to fix our notation
and state the established results. The second one contains all
the proofs corresponding to the analysis of the state system, while the last one 
is completely devoted to the control problem. Namely, the last section faces the analysis of
the existence of optimal control,
the linearized problem, the investigation of the \Frechet~differentiability of the control-to-state
mapping and the adjoint problem. Moreover, it
contains the necessary conditions that a control has to satisfy to be optimal.

\section{General assumptions and results}
\label{SEC_RESULT}
\setcounter{equation}{0}
In the following, we intend to fix the notation, state the problem in a 
precise form, and announce the main results.

The introduction should not have created any confusion since
the employed notation is quite standard. 
We assume $\Omega$ to be a smooth, bounded and connected open 
set in $\erre^3$, whose boundary is denoted by $\Gamma$. From the smoothness property, it
is almost everywhere well defined the unit normal vector $n$ of $\Gamma$ and the symbol $\dn$ represents
the outward derivative in that direction.
Moreover, for a fixed $T > 0$, which stands for the final time involved in the 
evolution process, we set
\Bsist
	& \non
	Q_{t}:=\Omega \times (0,t) \quad \hbox{and} \quad \Sigma_{t}:=\Gamma \times (0,t)
	\quad \hbox{for every $t \in (0,T]$,}
	\\ \non
	& Q:=Q_{T}, \quad \hbox{and} \quad \Sigma:=\Sigma_{T}.
\Esist
As the functional spaces are concerned, it turns out to be 
very convenient to introduce the following
\Beq
	\non
	H:= \Ldue, \quad V:= \Huno, \quad W:=\graffe{v \in \Hdue : \dn v = 0 \hbox{ on } \Gamma},
\Eeq
and endow them with their standard norms indicated by $\norma{\cdot}_{\bullet}$, 
where $_\bullet$ stands for the referred space or is completely omitted if it is clear from the 
context which norm should be. In the same way, we write $\norma{\cdot}_{p}$ for the usual norm 
in $L^p(\Omega)$. The above definitions yield that $(V,H,\Vp)$ forms a Hilbert triplet,
that is, the following injections $V \subset H \equiv H^* \subset V^*$ 
are both continuous and dense. As a consequence, we also have that
$ \< u,v > = \iO uv$ for every $u \in H$ and $v \in V$, where $\< \cdot , \cdot >$ 
denotes the duality pairing between the dual $\Vp$ and $V$ itself.

Now, we state the general assumptions on the problem.
\begin{enumerate}[label=\bf{(H\arabic*)}]
	\item 	$\bz, \bQ, \bO,\bQh,\bOh $ are nonnegative constants, but not all zero.
	\label{H1}
	\item 	$\phQ, \andrea{\sQ} \in L^2(Q), \phO, \andrea{\sO} \in \Huno,$
			$u_*,u^*\in L^\infty(Q) \hbox{ with } u_*\leq u^* \,\aeQ.$
	\label{H2}
	\item   $\a, \b > 0$.
	\label{H3}	
	\item 	$\m_0 \in \Huno \cap L^\infty(\Omega) , \ph_0 \in \Hdue, \s_0 \in \Huno$.
	\label{H4}
	\item 	$P \in C^2(\erre)$ is nonnegative, bounded and Lipschitz continuous. 
	\label{H5}
	\item   $\hat{B}:\erre \to [0, \infty]$ is convex, proper and lower semicontinuous, with $\hat{B}(0)=0$.
	\label{H6}
	\item   $\hat{\pi} \in C^3(\erre)$ and $\pi:=\hat{\pi}'$ is 		
			\Lip~continuous.
	\label{H7}
\end{enumerate}

We define the potential $ F: \erre \to [0,\infty]$ and the graph 
$B \subseteq \erre \times \erre$ by 
\Beq
\label{Fpot}
F:= \hat{B} + \hat{\pi} \quad \hbox{and} \quad B:=\partial\hat{B},
\Eeq
and note that $B$ is a maximal monotone operator (see, e.g., \cit{\cite[Ex. 2.3.4, p. 25]{BRZ}})
with domain denoted by $D(B)$.
Furthermore, we assume that $B$, when restricted to its domain $D(B)$, is a smooth function.
Indeed, we require that
\begin{enumerate}[label=\bf{(H\arabic*)}]
	\setcounter{enumi}{7}
	\item  $\ D(B)=(r_-,r_+), \hbox{ with } -\infty \leq r_- <0 < r_+ \leq + \infty,$ 
	$B(0)=0$, 
	\\ $F_{| _{D(B)}} \in C^3(r_-,r_+)$, 
	\hbox{ and } $\lim\limits_{r \to r_{\pm}} F'(r) = \pm \infty$.
	\label{H8}
	\item 	$r_- < \inf \ph_0 \leq \sup \ph_0 < r_+$.
	\label{H9}
	\item   $1/\b \tonde{\m_0 + \Delta \ph_0 - B(\ph_0) - \pi(\ph_0)} \in \Ldue.$
	\label{H10}
\end{enumerate}
It is worth to underline that from the above requirements, it follows that
both $\hat{B}(\ph_0)$ and $B(\ph_0)$ are both in $\Linfty$, thus a fortiori in $\Luno$.
In the literature, with a slight abuse of notation, $F'$ usually denotes
the sum of $B$, the subdifferential of $\hat{B}$, and $\pi$, namely $F'=B+\pi$.
Here, since $F$ is regular, $B$ exactly represents the derivative of $\hat{B}$
in $(\rmin,\rmax)$.

Notwithstanding \ref{H6}-\ref{H10}, let us point out that there are
significant classes of double-well potentials that fit the assumptions.
Standard choices are the regular potential and the, physically more relevant, logarithmic 
one. Written as \eqref{Fpot}, they read as
\Bsist
	& F_{reg}(r):= \frac 14 (r^2-1)^2 
	= \frac 14 r^4 - \frac 14 (2 r^2-1)
	\quad \hbox{for } r \in \erre,
	\label{regpot}
	\\ 
	& F_{log}(r):= 
	((1-r)\log(1-r)
	+(1+r)\log(1+r))
	- k r^2
	\quad \hbox{for } |r|<1,
	\label{logpot}
\Esist
where in the latter $k$ is a constant large enough to kill convexity.
Moreover, it is usually helpful to extend \eqref{logpot} by continuity
imposing that it assumes the value $+\infty$ outside its actual domain.
Note that both \eqref{regpot} and \eqref{logpot} do fit our framework
with $D(B)=(-\infty,+\infty)$ and $D(B)=(-1,+1)$, respectively.
Furthermore, if we take into account $F_{reg}$, due to its regularity,
all the results we are going to prove still hold true even in a slightly weaker framework.
However, since $F_{reg}$ is introduced as an approximation of more general potentials,
we try to focus our attention on the singular ones, such as
$F_{log}$, which is more relevant for the applications.
Before starting with the statements, we introduce another notation.
\Beq
	\non
	\hbox{Let } \UR \hbox{ be an open set in } L^2(Q) \hbox{ such that }
	\Uad \subset \UR \hbox{ and } \norma u_2 \leq R \hbox{ for all }
	u \in \UR.
\Eeq
As it usually occurs in control problems, the requirements \ref{H2}-\ref{H10} 
are far from sharp in terms of the well-posedness and regularity result 
of \EQ~are concerned. Anyhow, they are all useful in order to deal with the 
corresponding control problem. 

Let us proceed this section by listing the obtained results.
\Bthm[\bf well-posedness and separation results]
\label{THMstate}
Under the hypotheses \ref{H2}-\ref{H10}, and for every $u \in \UR$, 
the following results hold true.\\
$(i) \,\,$ The system \EQ~has a unique strong solution $(\m,\ph,\s)$
which satisfies
\Bsist
	&\ph \in W^{1,\infty}(0,T;H) \cap \H1 V \cap \L\infty W \subset \C0 {C^0({\bar{\Omega}})}
	\label{reg1}
	\\ 
	&\m, \s \in \H1 H  \cap \L\infty V \cap \L2 W \subset C^0([0,T];V)
	\label{reg2}
	\\ 
	&\m \in L^{\infty}(Q)
	\label{reg3}
\Esist
\Accorpa \reg reg1 reg3
that is, there exists a constant $C_1>0$, which depends on $R$, $\a$ and $\b$, and on the data 
of the system, such that
\Bsist
 	&& \non 
 	\norma{\ph}_{W^{1,\infty}(0,T;H) \cap \H1 V \cap \L\infty W}
	+ \norma{\m}_{\H1 H  \cap \L\infty V \cap \L2 W \cap L^\infty(Q)}
	\\ && \quad
	+ \norma{\s}_{\H1 H \cap \L\infty V \cap \L2 W}
	\leq C_1.
	\label{regstima}
\Esist
$(ii) \,\,$ There exists a compact subset $K$ of $(r_-,r_+)$ such that
\Beq
	\ph(x, t) \in K \quad \hbox{for all } (x,t) \in Q;
\Eeq
in particular, there exists a constant $C_2>0$, which depends on $R$, $\a$ and $\b$, 
$K$ and on the data of the system, such that
\Beq
	\label{stimasep}
	\norma{\ph}_{C^0(\overline{Q})}
	+ \max_{0 \leq i \leq 3} \,\norma{F^{(i)}(\ph)}_{L^\infty(Q)}
	+ \max_{0 \leq j \leq 2} \,\norma{P^{(j)}(\ph)}_{L^\infty(Q)}
	\leq C_2.
\Eeq
\Ethm

\Bthm[\bf continuous dependence on the control]
\label{THcontdep}
Assume \ref{H2}-\ref{H10}. Then there exists a constant $C_3 > 0 $, which depends
only on $R$, $\a$ and $\b$, and on the data of the system such that, if $u_i \in \UR$
and $(\m_i,\ph_i,\s_i)$ are the corresponding solutions with the same initial value, 
$i=1,2$, it holds
\Bsist
	\non &&
	\norma{\a(\m_1-\m_2) + (\ph_1-\ph_2) + (\s_1-\s_2)}_{\L\infty\Vp}
	+ \norma{\m_1-\m_2}_{\L2H}
	\\ \non && \quad
	+\norma{\ph_1-\ph_2}_{\L\infty H \cap \L2 V}
	+ \, \norma {\s_1-\s_2} _{\L\infty H \cap \L2 V}
	\\ && \quad
	\leq C_3 \norma {u_1-u_2} _{\L2 H}.
	\label{contdepcontrol}
\Esist
\Ethm
Let us remark that in the proof of the above result we do not account for
point $(ii)$ of Theorem \ref{THMstate}.
In fact, we will see that this first continuous dependence result is not sufficient to handle
the \CP~(particularly to prove the \Frechet~differentiability of 
the control-to-state mapping $\S$, cf. Sec. \ref{FRECHET}), then in the beneath 
lines there is an improvement that, this time, take strongly into account 
the second part of Theorem \ref{THMstate}.

\Bthm
\label{THcontdep2}
In the same framework of Theorem \ref{THcontdep}, there exists a constant 
$C_4 > 0 $, possibly smaller than $C_3$, which depends
only on $R$, $\a$ and $\b$, and on the data of the system such that 
\Bsist
	\non &&
	\norma{\m_1-\m_2}_{\L\infty H \cap \L2 V}
	+ \norma{\ph_1-\ph_2}_{\H1 H \cap \L\infty V}
	\\ && \quad
	\leq C_4 \norma {u_1-u_2} _{\L2 H}.
	\label{contdepcontroldue}
\Esist
\Ethm

Since we have already provided the well-posedness of the state system in Theorem \ref{THMstate}, 
we can introduce the \socal~control-to-state mapping $\S$ that will cover a central
role in the control theory. It consists of the map that assigns to every
admissible control $u$ the corresponding solution triple $(\m,\ph,\s)$, which components 
belong to the functional spaces pointed out by \reg. 
Moreover, it allows us to present the \socal~reduced cost functional as follows
\Bsist
	\non &
	\tilde{\J} : \UR \to \erre, \quad \hbox{defined by} \quad
	\tilde{\J} (u) := \J({\S}_{2,3}(u),u), 
	\\ &   \non
	\hbox{where } {\S}_{2,3}(u) 
	\hbox{ represents the couple of the second 
	and third components} 
	\\ & \hbox{of the solution triple ${\S}(u)=(\m,\ph,\s)$.}
	\label{reducedcontrol}
\Esist
In this view, Theorem \ref{THcontdep} established the 
\Lip~continuity of $\S$ in this natural functional framework.

At this point, we introduce a well-posedness result for the linearized system, 
which comes out naturally from the investigation of the control problem. 
First of all, let us present the mentioned problem. 
Fixed $\bar{u} \in \UR$, we denote $(\bm,\bph,\bs)=
{\S}(\bar{u})$ the corresponding solution to \EQ.
Then, for any $h \in L^2(Q)$, the linearized system reads as
\Bsist
  & \a \dt \et + \dt \th - \Delta \et = P'(\bph) (\bs - \bm)\th + P(\bph)(\r - \et)
  \quad \hbox{in $\, Q$}
  \label{EQLinprima}
  \\[0.2cm]
  & \et = \b \dt \th - \Delta \th + F''(\bph)\th
  \label{EQLinseconda}
  \quad \hbox{in $\,Q$}
  \\[0.2cm]
  & \dt \r - \Delta \r = -P'(\bph) (\bs - \bm)\th - P(\bph)(\r - \et) + h
  \label{EQLinterza}
  \quad \hbox{in $\,Q$}
  \\[0.3cm]
  & \dn \r = \dn \th = \dn \et = 0
  \quad \hbox{on $\,\Sigma$}
  \label{BCEQLin}
  \\[0.2cm]
  & \r(0) = \th(0) = \et (0) = 0
  \quad \hbox{in $\,\Omega$}.
  \label{ICEQLin}
\Esist
\Accorpa\EQLin EQLinprima ICEQLin
Here the existence and uniqueness result follows.
\Bthm[\bf well-posedness of the linearized system]
\label{THMlinear}
Under the assumptions \ref{H2}-\ref{H10}, and for every h $\in L^2(Q)$,
the system \EQLin~possesses a unique solution triple $(\et,\th,\r)$
which satisfies
\Bsist
	\et, \th, \r \in \H1 H \cap \L\infty V \cap \L2 W \subset \C0 V;
	\label{reglin}
\Esist
that is, there exists a constant $C_5>0$, which depends on the data 
of the system, and possibly on $\a$ and $\b$, such that
\Bsist
 	\non &&
 	\norma{\et}_{\H1 H \cap \L\infty V \cap \L2 W}
	+ \norma{\th}_{\H1 H \cap \L\infty V \cap \L2 W}
	\\ && \non \quad
	+ \norma{\r}_{\H1 H \cap \L\infty V \cap \L2 W}
	\leq C_5.
\Esist
\Ethm

In the following, we prove that 
$\S$ is even \Frechet~differentiable in suitable Banach spaces.
\Bthm[\bf \Frechet~differentiability of $\S$]
\label{THMFrechet}
Assume \ref{H2}-\ref{H10}. Then the control-to-state mapping $\S$ is 
\Frechet~differentiable in $\UR$ as a mapping from $L^2(Q)$ into the state space
$\cal Y$, where 
\Beq
	\label{statespaceFre}
	{\cal Y}:= \biggl( \H1 H \cap \L\infty V \cap \L2 W \biggr)^3.
\Eeq
Moreover, for any $\bar{u} \in \UR$, the \Frechet~derivative 
$D\S(\bar{u})$ is a linear and continuous operator from $L^2(Q)$ to $\calY$,
and for every $ h \in L^2(Q)$,
$D\S(\bar{u}) h = (\et, \th, \r)$ where $(\et, \th, \r)$ is the unique solution
to the linearized system \EQLin~associated with $h$.
\Ethm

\Bthm[\bf Existence of optimal control]
\label{THexistenceofcontrol}
Assume \ref{H2}-\ref{H10}. Then the optimal control problem \CP~has at least a solution
$\bar{u} \in \Uad$.
\Ethm

As the necessary optimality condition is concerned,
we recall the reduced cost functional \eqref{reducedcontrol} and the fact that $\Uad$ is convex. 
Therefore, the optimal inequality we are looking for turns out to be
\Beq
	\label{abstrnec}
	\< D\tilde{\J}(\bar{u}) , v - \bar{u} > \geq 0 \quad \hbox{for every } v \in \Uad,
\Eeq 
where $ D\tilde{\J}(\bar{u})$ represents the differential of $\tilde{\J}$, at least in the G\^ateaux sense.
Accounting for Theorem \ref{THMFrechet} and the chain rule, \eqref{abstrnec} develops as follows.

\Bcor
\label{CORprimanec}
Suppose that the assumptions \ref{H1}-\ref{H10} are fulfilled.
Let $\bar{u} \in \Uad$ be an optimal control for \CP~with his corresponding optimal state
$(\bm,\bph,\bs)=\S(\bar{u})$. Then we have
\Bsist
  \non
  && \bQ \intQ (\bph - \phQ)\th 
  + \bO \iO (\bph(T) - \phO)\th(T)
  + \andrea{\bQh \intQ (\bs - \sQ)\r }
  \\ && \quad
  + \andrea{\bOh \iO (\bs(T) - \sO)\r(T)}
  + \bz \intQ \bar{u}(v - \bar{u}) 
  \geq 0 \quad \forall v \in \Uad,
    \label{primanec}
\Esist
where $\th$ and $\r$ are the second and third components of the unique solution triple $(\et,\th,\r)$ to the linearized 
system \EQLin~associated with $ h = v -\bar{u}$.
\Ecor

To eliminate the presence of the variables $\th$ and $\r$ in the previous inequality,
we introduce the \socal~adjoint problem that consists of the following 
system of partial differential equations.
\Bsist
  & \b \dt q - \dt p + \Delta q - F''(\bph)q + P'(\bph)(\bs - \bm)(r-p)= \bQ(\bph - \phQ)
  \quad \hbox{in $\, Q$}
  \label{EQAggprima}
  \\[0.2cm]
  & q -\a \dt p - \Delta p + P(\bph)(p -r) = 0
  \label{EQAggseconda}
  \quad \hbox{in $\,Q$}
  \\[0.2cm]
  & -\dt r - \Delta r + P(\bph)(r - p) = \andrea{\bQh (\bs - \sQ)}
  \label{EQAggterza}
  \quad \hbox{in $\,Q$}
  \\[0.4cm]
  &\dn q = \dn p = \dn r = 0
  \quad \hbox{on $\Sigma$}
  \label{BCEQAgg}
  \\[0.2cm] \non
  & p(T) - \b q(T) =\bO(\bph(T) - \phO),  \quad
  \\ & \quad
  \a p(T) = 0, 
  \quad r(T) = \andrea{\bOh (\bs(T) - \sO)}
  \quad \hbox{in $\Omega$.}
  \label{ICEQAgg}
\Esist
\Accorpa\EQAgg EQAggprima ICEQAgg
Here the existence and uniqueness result for the adjoint problem is stated.

\Bthm[\bf Well-posedness of the adjoint problem]
\label{THMadjoint}
Under the assumptions \ref{H1}-\ref{H10}, the system \EQAgg~has a unique 
solution $(q,p,r)$ that satisfies the following regularity requirements
\Bsist
	q,p,r \in \H1 H \cap \L\infty V \cap \L2 W \andrea{\subset \C0 V}.
	\label{regadj}
\Esist
\Ethm

Finally, the well-posedness of the adjoint system allows us to improve Corollary 
\ref{CORprimanec} leading
to a second necessary condition. Namely, we achieve the following result.

\Bthm[\bf Necessary optimality condition]
\label{THMsecondanec}
Assume \ref{H1}-\ref{H10}. Let $\bar{u} \in \Uad$ be an optimal control 
with his corresponding optimal state $(\bm,\bph,\bs)=\S(\bar{u})$
and let $(p,q,r)$ be the solution to the corresponding adjoint system. Then we have
\Bsist
  \label{secondanec}
  \intQ (r + \bz \bar{u})(v - \bar{u})
  \geq 0 \quad \forall v \in \Uad.
\Esist
\Ethm

To conclude the section, let us introduce further notation and recall
some \wk~inequalities and general facts related to the \CH~equation.
First of all, we remind the Young inequality
\Beq
  ab \leq \delta a^2 + \frac 1 {4\delta} \, b^2
  \quad \hbox{for every $a,b\geq 0$ and $\delta>0$}.
  \label{young}
\Eeq
Furthermore, for given $v\in\Vp$ and $\underline v\in\L1\Vp$,
we introduce their \generaliz ed mean values 
$\vO\in\erre$ and $\underline v^\Omega\in L^1(0,T)$ by
\Beq
  \vO := \frac 1 {|\Omega|} \, \< v , 1 >,
  \aand
  \underline v^\Omega(t) := \bigl( \underline v(t) \bigr)^\Omega
  \quad \aat,
  \label{media}
\Eeq
where \eqref{media} reduces to the usual mean values when it is applied to elements
of~$H$ or $\L1H$.
In addition, we often owe to the 
Poincar\'e inequality
\Beq
  \normaV v^2 \leq C_{\Omega} \bigl( \normaH{\nabla v}^2 + |\vO|^2 \bigr)
  \quad \hbox{for every $v\in V$},
  \label{poincare}
\Eeq
where we stressed the fact that $C_\Omega$ depends on $\Omega$. 
Since it will be convenient to interpret some partial differential equations
in the framework of the Hilbert triplet $(V,H,\Vp)$, we introduce the Riesz isomorphism 
associated with $V$. That is, we define the map $\A:V \to \Vp$ as follows
\Beq
	\< \A u,v> = (u , v)_{V} = \iO (\nabla u \cdot \nabla v + u v)  \quad \hbox{for every } u, v \in V.
	\label{defA}
\Eeq
Observe that when restricted to its domain $W$, $\A$ turns out to be the 
operator $- \Delta + I$ endowed with homogeneous Neumann boundary conditions, where $I$ denotes
the identity map of $W$. A little investigation on $\A$ leads to the following identities
\Bsist
	&& \<\A u, \A^{-1} v^*> = \< v^*, u > \quad \hbox{for all $u \in V$ and $v^* \in \Vp$,}
	\label{propAuno}
	\\ 
	&& \<u^*, \A^{-1} v^*> = ( u^*, v^*)_* \quad \hbox{for all $u^*, v^* \in \Vp$,}
	\label{propAdue}
\Esist
where $(\cdot,\cdot)_*$ stands for the inner product of $\Vp$,
whence also
\Beq
  2 \< \dt v^*(t) , \A^{-1} v^*(t) >
  = \frac d{dt} \, \normaVp{v^*(t)}^2
  \quad \aat.
  \label{dtA}
\Eeq
\Accorpa\propA defA dtA

\Brem
\label{convcostanti}
Let us explain a convention that we are going to use throughout the paper. Since we have to 
deal with a lot of estimates, we agree that the symbol $c$ stands for any constants 
which depend only on the final time~$T$, on~$\Omega$, the shape of the nonlinearities,
on the norms of the involved functions, and possibly on 
$\a$ and $\b$. For this reason, the meaning of $c$ might
change from line to line and even in the same chain of inequalities.
Conversely, the capital letters are devoted to denote precise constants.
\Erem

\section{State system and continuous dependence results}
\label{SEC_STATESTYSTEM}
\setcounter{equation}{0}
From this section on, we will focus our attention to prove the statements.
This section is devoted to the investigation of the state system, namely
we aim at checking Theorems \ref{THMstate}, \ref{THcontdep}, and \ref{THcontdep2}. 
Let us begin dealing with the first one.

\proof[Proof of Theorem \ref{THMstate}]
In \cit{\cite[Thm. 2.2, p. 2426]{CGH}} it has been shown 
that the system \EQ~possesses a unique strong solution
with the following regularity
\Bsist
	\non
	\m,\ph,\s \in \H1 H \cap \L2 W \subset \C0 V,
\Esist
in the homogeneous case $u \equiv 0$.
Since we admit that $u$ can be chosen in $\UR$, only \sfw~modifications
are needed in order to prove that, for every choice of $u$ in $\UR$,
there exists a unique corresponding solution $(\m,\ph,\s)$ satisfying the same 
regularity mentioned above.
Let us point out that conditions \ref{H2}-\ref{H10} perfectly fit the framework of \cite{CGH}.
In fact, the strong requirement $(2.6)$ of \cite{CGH} is only needed to 
handle the asymptotic behavior, whereas it can be substituted by the weak requirement \ref{H7} as
the investigation of the well-posedness and regularity of the system are concerned.

In the following, we proceed formally; as a matter of fact, we should introduce suitable 
approximation of the potential depending on a small parameter $\eps$ and then, 
after showing sufficient compactness property, let $\eps\searrow0$
as made in \cite{CGH}. Anyhow, we will
take care in referring to works in which this strategy is properly employed to 
justify all the passages that we present only in a formal level.

With the following estimates, we aim at improving the regularity of the unique solution 
to \EQ~in view of the forthcoming control investigation. 
Once obtained, it is a standard matter 
to conclude by compactness arguments that the solution triple satisfies 
\accorpa{reg1}{reg3}. 
The rigorous treatment of the first three estimates, with little variations, can be found in 
\cit{\cite[eqs. (4.4)-(4.12), pp. 2431-2432]{CGH}}.

{\bf First estimate:} 
We add to both the sides of \eqref{EQseconda} the term $\ph$,
multiply \eqref{EQprima} by $\m$, this new second equation by $-\dt\ph$ and \eqref{EQterza}
by $\s$, then we add the resulting equations and integrate over $Q_t$ and by parts.
A little rearrangements of the terms produce
\Bsist
	\non &&
	\frac \a2 \IO2 \m
	+ \I2 {\nabla\m}
	+ \b \I2 {\dt\ph}
	+ \frac 12 \IO2 {\ph}
	+ \frac 12 \IO2 {\nabla\ph}
	\\ \non && \quad
	+ \iO \hat{B}(\ph(t))
	+ \frac 12 \IO2 \s
	+ \I2 {\nabla\s}
	+ \intQt P(\ph)(\s-\m)^2
	\\ \non && \quad
	= 
	\frac \a2 \iO |\m_0|^2
	+ \frac 12 \iO |\ph_0|^2
	+ \frac 12 \iO |\nabla\ph_0|^2
	+ \iO \hat{B}(\ph_0)
	\\ \non && \quad
	+ \frac 12 \iO |\s_0|^2
	+ \intQt u\s
	+ \intQt \ph\,\dt\ph
	- \intQt \pi(\ph) \dt \ph,
\Esist
where we split $F'$ as sum of $B$ and $\pi$, and where the former, 
multiplied by $\dt\ph$, consists of the derivative with respect to time of $\hat{B}(\ph(t))$.
All the terms on the \lhs~are nonnegative since they are squares and  
$P$ and $\hat{B}$ are nonnegative by \ref{H5} and \ref{H6}, respectively.
The first five terms on the \rhs~are easily managed owing to \ref{H4}
and to the properties of $\hat{B}$, while the others were denoted by $I_1,I_2,I_3$.
Accounting for the Young inequality \eqref{young}, we obtain
\Bsist
	\non
	|I_1|+|I_2|+|I_3|
	\leq
	\frac 12 \I2 u
	+ \frac 12 \I2 \s
	+ 2\d \I2 {\dt\ph}
	+ \cd \intQt  (|\phi|^2+1).
\Esist
We choose $0< \d < \b/2$, and invoke the Gronwall lemma to conclude that
\Bsist
 	\non &&
 	\norma{\m}_{\L\infty H \cap \L2 V}
	+ \norma{\ph}_{\H1 H \cap \L\infty V}
	+ \norma{\hat{B}(\ph)}_{\L\infty \Luno}
	\\ \quad &&
	+ \norma{\s}_{\L\infty H \cap \L2 V}
	\leq 
	c.
	\label{Iest}
\Esist

{\bf Second estimate:} 
Now we multiply \eqref{EQprima} by $\dt\m $ and \eqref{EQterza}
by $\dt\s$, add the resulting equations and integrate over $Q_t$.
Using \eqref{Iest} and the boundedness of $P$, and arguing as above
lead to
\Bsist
	&& \norma{\m}_{\H1 H \cap \L\infty V}
	+ \norma{\s}_{\H1 H \cap \L\infty V}  
	\leq
	c.
	\label{IIest}
\Esist

{\bf Third estimate:} 
Equations \eqref{EQprima} and \eqref{EQterza} show a parabolic structure with respect to
$\m$ and $\s$, respectively. Moreover, it follows from the previous estimates
that their forcing terms are both in $\L2 H$.
Therefore, since the initial data \eqref{ICEQ} are in $V$ (cf. \ref{H4}),
parabolic regularity theory with Neumann homogeneous boundary conditions gives
\Bsist
	\norma{\m}_{\L2 W}
	+ \norma{\s}_{\L2 W}
	\leq 
	c.
	\label{IIIest}
\Esist

{\bf Fourth estimate:} 
As above we would like to obtain more regularity for the phase variable $\ph$ by comparing terms in
\eqref{EQseconda}. Since it is more delicate, worth to show the detailed procedure.
In fact, we can rearrange \eqref{EQseconda} as follows
\Bsist
	\label{fortheq}
	-\Delta \ph + B(\ph) = f, \quad \hbox{where $f:=\m -\b \,\dt \ph -\pi(\ph)$.}
\Esist
The previous estimates entails that $f \in \L2 H$.
Now, we multiply the above inequality by $(- \Delta \ph)$ and integrate over $\Omega$.
Actually, note that this choice is rigorously forbidden since too few regularity 
is known on the phase variable. Anyhow, this choice can be formally justified 
by introducing a suitable Faedo-Galerkin scheme.
So, \aat, we get the following inequality
\Bsist
	\non
	\IO2 {\Delta\ph}
	+\iO B'(\ph(t)) \, |\nabla \ph(t)|^2 
	\leq 
	- \iO f(t)\,\Delta\ph(t)
	\leq 
	\frac 12 \IO2 {\Delta\ph}
	+\frac 12 \IO2 {f},
\Esist
owing to \eqref{young}.
Both the terms on the \lhs~are nonnegative since $B'$ is so.
Hence, we realize that 
\Beq
	\non
	\norma{\Delta\ph}_{\L2 H} \leq c.
\Eeq
Moreover, by elliptic regularity theory, the boundary conditions \eqref{BCEQ}, 
and by comparison in \eqref{fortheq}, we conclude that
\Beq
	\norma{\ph}_{\L2 W} + \norma {B(\ph)}_{\L2 H}
	\leq c.
	\label{IVest}
\Eeq

{\bf Fifth estimate:} 
We continue to proceed formally, in order to keep the proof as short and easy as possible.
Here, for a precise and detailed treatment it will be necessary to introduce time steps and suitable 
translations, and show some estimates for this new functions. This 
procedure will become quite technical. Anyhow, for the interested reader,
we refer to \cit{\cite[Proof of Thm. 2.6 (iii), p. 2436]{CGH}}, where the correct 
procedure is performed to establish a slightly different estimate.

So, we differentiate \eqref{EQseconda} with respect to the time variable,
multiply it by $\dt\ph$, and integrate over $Q_t$ to get
\Bsist
	\non
	\intQt \dt\m \ \dt\ph 
	= \b \intQt \dtt\ph \ \dt\phi
	- \intQt (\Delta \dt\ph) \, \dt\ph
	+ \intQt (B'(\ph)+\pi'(\ph)) \ |\dt\ph|^2.
\Esist
Using the integration by parts and the boundary conditions 
\eqref{BCEQ}, we deduce that 
\Bsist
	&& \non
	\frac \b 2 \iO |\dt\ph(t)|^2
	+ \I2 {\nabla\dt\ph}
	+ \intQt B'(\ph)|\dt\ph|^2
	=
	\frac \b 2 \iO |(\dt\ph)(0)|^2
	\\ \non && \quad
	-\intQt \pi'(\ph)|\dt \ph|^2
	+ \intQt \dt\m \ \dt\ph,
\Esist
where the terms on the \lhs~are all nonnegative.
The first term of the \rhs~is under control, due to \eqref{EQseconda} and \ref{H10}.
Moreover, the last two integrals can be estimate as follows 
\Bsist
	\non
	-\intQt \pi'(\ph)|\dt \ph|^2
	+ \intQt \dt\m \ \dt\ph 
	\leq 
	c \intQt |\dt\ph|^2
	+\frac 12 \intQt |\dt\m|^2,
\Esist
owing to the \Lip~continuity of $\pi'$ and \eqref{young}.
Thus, thanks to \eqref{Iest} and \eqref{IIest} we obtain
\Bsist
	\norma{\ph}_{W^{1,\infty} (0,T;H) \cap \H1 V}
	\leq 
	c.
	\label{Vest}
\Esist

{\bf Sixth estimate:} 
We take again into account equation \eqref{fortheq}. Due to the previous
estimates, we can infer that $f$ is more regular than we pointed out before.
In fact, now we have that $f\in \L\infty H$. Therefore, the test by $-\Delta \ph$
leads to the estimate $\norma{\Delta\ph}_{L^{\infty}(0,T;H)} \leq c $. 
Moreover, by the boundary conditions, the elliptic regularity and 
comparison in \eqref{fortheq}, we deduce that
\Beq
	\norma{\ph}_{\L\infty W} 
	+ \norma{B(\ph)}_{\L\infty H}
	\leq c,
	\label{VIest}
\Eeq
which gives, by the Sobolev embeddings, also
\Beq
	\label{phinf}
	\norma{\ph}_{L^{\infty} (Q)} 
	\leq c.
\Eeq

Furthermore, an application of the \wk~embedding results (see e.g., \cit{\cite[Sect. 8, Cor. 4]{Simon}})
directly recovers the continuity of the solution variables. 
Namely, as the variables $\m$ and $\s$ are concerned,
due to \accorpa{Iest}{IIIest}, we infer that they belong to $\C0 V$.
Since the phase variable $\ph$ satisfies, in addition, the estimates
\eqref{IVest} and \accorpa{Vest}{phinf}, we deduce that $\ph$ is more 
regular and it belongs to $\C0 {C^0(\bar{\Omega})}$.

Now, we start to approach the separation result $(ii)$. This property will be 
crucial in order to handle the potential and its higher order derivatives.
Indeed, if $(ii)$ holds true, it acts on functions
which values are well detached from the boundary of the domain of $B$. 
In this way $F$ and his higher order derivatives do not blow up and
they turn out to be \Lip~continuous~and bounded functions.

First of all, we need to show the boundedness of the chemical potential in the whole of $Q$.
In this direction, we would like to apply \cit{\cite[Thm. 7.1, p.\ 181]{LDY}}
to equation \eqref{EQprima}. 
The key point is the parabolic structure with respect to $\m$ that \eqref{EQprima} possesses.
Indeed, by simply rearranging the terms, we get
\Beq
	\non
	\a \dt \m  - \Delta \m = g, \quad \hbox{where $g := P(\phi) (\sigma - \mu)	- \dt \ph $.}
\Eeq
Roughly speaking, the result formalizes the following idea: if the initial
datum is bounded in $\Omega$ and the forcing term $g$ satisfies a suitable summability
regularity with respect to space and time, then it is natural to expect that the variable $\m$ stay 
bounded in the whole of $Q$.
Actually, from \ref{H4} the property on the initial data is already satisfied.
Moreover, \ref{H5} and the previous estimates immediately yield that
\Beq
	\non
	\norma {g}_{\L\infty H} \leq c.
\Eeq
This allows us to apply \cit{\cite[Thm. 7.1, p.\ 181]{LDY}} and infer
that there exists a positive constant $c$ such that
\Beq
	\label{muinf}
	\norma{\m}_{L^{\infty}(Q)} \leq c.
\Eeq
Note that \eqref{muinf} ends the proof of $(i)$ and turns out to be fundamental in order to 
proceed with the second part of Theorem \ref{THMstate}.

As before, let us emphasize that the estimate
we are going to prove in the following is formal.
Anyhow, it can be reproduced correctly by introducing a suitable approximation
scheme, as made in \cit{\cite[Proof of Thm. 2.6, pp.~992-994]{CGS}}.

{\bf Seventh estimate:} 
We multiply \eqref{EQseconda} by 
$|B(\ph)|^{p-1}\, \sign \ph = |B(\ph)|^{p-2} B(\ph)$, for a fixed $p>2$, 
and integrate over $Q_t$. 
Moreover, we set $f:= \m - \pi(\ph)$ and observe that $f$ belongs to $L^\infty(Q)$
due to \eqref{phinf}, \eqref{muinf} and \ref{H7}.
We infer, for every $t\in[0,T]$, that
\Bsist
	&& \non
	\b \iO {\cal B}_p (\ph(t))	
	+ (p-1) \intQt B'(\ph)|B(\ph)|^{p-2}|\nabla \ph|^2
	+ \intQt |B(\ph)|^p 
	\\ && \quad
	= \b \iO {\cal B}_p (\ph_0)
	+\intQt f \, |B(\ph)|^{p-1} \sign \ph,
	\label{seventhest}
\Esist
where ${\cal B}_p(r):= \int _0^r |B(s)|^{p-1}\sign s \,ds$.
\dafare{Furthermore, all the terms on the \lhs~are nonnegative.}
As the \rhs~is concerned, we manage the first term in the following way.
From \ref{H9} we know that $|B(\ph_0)|$ is bounded by a positive constant $M$,
hence, we infer that
\Beq
	\non
	\b \iO {\cal B}_p (\ph_0) \leq 
	\b M^{p-1} \dafare \iO{|\ph_0|}\leq c^p.
\Eeq
Moreover, the last term can be estimated by
\Bsist
	\non &&
	\intQt f \, |B(\ph)|^{p-1} \sign \ph
	\leq \frac 1p\, c^p
	+ \frac 1{p'} \intQt |B(\ph)|^{(p-1){p'}}
	\leq
	c^p
	+ \frac 1{p'} \intQt |B(\ph)|^p,
\Esist
owing to the general version of the Young inequality, 
where $p'$ stands for the conjugate exponent of $p$.
Using the above estimates, we can rearrange \eqref{seventhest} to conclude that
\Beq
	\non
	\frac 1p \intQt |B(\ph)|^p \leq c^p,
\Eeq
that implies 
\Beq
 	\non
 	\norma{B(\ph)}_{L^p(Q)}\leq c,
\Eeq
where the constant $c$ is independent of $p$. Since the above procedure can
be iterated for every $p>2$, we realize that $\norma{B(\ph)}_{L^\infty(Q)}\leq c$
and from this we recover $\norma{F'(\ph)}_{L^\infty(Q)}\leq c$.
In view of \ref{H8}, this establishes that
\Beq
	\non
	r_- < \inf \ph \leq \sup \ph < r_+ \quad \aaQ,
\Eeq
as we claimed. \qed 

At this point, we prove the continuous dependence results.
\proof[Proof of Theorem \ref{THcontdep}]
For this proof we have largely taken inspiration from \cit{\cite[Sec. 3, pp. 2429-2430]{CGH}}.
First of all, we set 
\Beq
	\label{notdiff}
	\m:= \m_1-\m_2, \quad \ph:= \ph_1-\ph_2, \quad \s:= \s_1-\s_2, \quad u:= u_1-u_2.
\Eeq
Writing \EQ~for $(\m_i,\ph_i,\s_i)$, $i=1,2$, and taking
the difference, we obtain the following equations and conditions:
\Bsist
  \label{FVdiprimaprova}
  & \a \dt \m + \dt \ph - \Delta \m = R_1-R_2
  \quad \hbox{in $\, Q$}
  \\
  \label{FVdifseconda}
  & \mu = \beta \dt \phi - \Delta \phi + F'(\phi_1) - F'(\phi_2)
  \quad \hbox{in $\,Q$}
  \\
  & \dt \sigma - \Delta \sigma = - (R_1-R_2) + u
  \quad \hbox{in $\,Q$}
  \label{FVdifterza}
  \\[0.1cm]
  &\dn \m = \dn \ph = \dn \s = 0
  \quad \hbox{on $\Sigma$}
  \\
  & \m(0) = \ph(0) = \s (0)= 0
  \quad \hbox{in $\Omega$}
\Esist
\Accorpa\EQContDep FVdiprimaprova FVdifterza
where $R_i:= P(\ph_i)(\s_i-\m_i)$, $i=1,2$.
Now, we take the sum of \eqref{FVdiprimaprova} and \eqref{FVdifterza}, 
then add to both the members of this new equation $\m+\s$. 
This gives
\Bsist
  & \dt ( \alpha \mu + \phi + \sigma) 
  + \A(\mu + \sigma)  
  = u + \mu + \sigma
  \quad \hbox{in $\, Q$},
  \label{FVdifprima}
\Esist
owing to \eqref{defA}.
Keeping in mind \propA, we multiply \eqref{FVdifprima} by $\A^{-1}(\alpha \mu + \phi +\sigma)$, 
\eqref{FVdifseconda} by $-\phi$ and \eqref{FVdifterza} by $\sigma$, add them, and integrate over $Q_t$.
We deduce that
\Bsist
	\non
	&& \frac 1 2 \normaVp {( \alpha \mu + \phi + \sigma)(t)}^2
	+ \intQt (\mu + \sigma) (\alpha \mu + \phi + \sigma)
	- \intQt \phi \mu
	\\ && \quad \non
	+ \frac \beta 2 \iO |\phi(t)|^2
	+ \intQt |\nabla\phi|^2
	+ \intQt (F'(\ph_1)-F'(\ph_2))\phi
	+ \frac 12 \iO |\sigma(t)|^2
	\\ && \quad \non
	+ \intQt |\nabla\sigma|^2
	= \iot \< u + \mu + \sigma, \A^{-1} (\alpha \mu + \phi + \sigma) >
	\\ && \quad
	- \intQt \tonde {P(\phi_1)-P(\phi_2)}(\sigma_1-\mu_1)\sigma
	- \intQt P(\phi_2) \tonde{\sigma - \mu} \sigma
	+ \intQt u \sigma,
	\label{depconstima}
\Esist
where the second and third terms of the \rhs~come from a simple rearrangement
of $R_1-R_2$.
We develop the second term of the \lhs~as 
\Bsist
	 \non
	 \alpha \intQt |\mu|^2 +\intQt |\sigma|^2 
	 + \intQt \mu \ph + (1+\alpha)\intQt \mu\sigma
	 + \intQt \sigma \phi,
\Esist
and move the last three terms of the above sum to the \rhs~of \eqref{depconstima}.
Observe that the term that involves the double-well potential should be decomposed as
\Bsist
	\non
	\intQt (F'(\ph_1)-F'(\ph_2))\ph 
	= \intQt (B(\ph_1) - B(\ph_2))\ph 
	+ \intQt (\pi(\ph_1) - \pi(\ph_2))\ph,
\Esist
where the first term of the \rhs~is nonnegative by the monotonicity of $B$, while the 
second one can be moved to the \rhs~of \eqref{depconstima} and easily managed, since $\pi$ is \Lip~continuous by \ref{H7}.
If we rearrange \eqref{depconstima} according to the above observations, we obtain
\Bsist
	\non 
	&& \frac 1 2 \normaVp {( \alpha \mu + \phi + \sigma)(t)}^2
	+ \alpha \intQt |\mu|^2 +\intQt |\sigma|^2 
	+ \frac \beta 2 \iO |\phi(t)|^2
	\\ && \quad \non	
	+ \intQt |\nabla\phi|^2
	+ \intQt (B(\ph_1) - B(\ph_2))\ph  
	+ \frac 12 \iO |\sigma(t)|^2
	+ \intQt |\nabla\sigma|^2
	\\ && \quad \non
	= \iot \< u + \mu + \sigma, \A^{-1} (\alpha \mu + \phi + \sigma) >
	- \intQt \tonde {P(\phi_1)-P(\phi_2)}(\sigma_1-\mu_1)\sigma
	\\ && \quad \non
	- \intQt P(\phi_2) \tonde{\sigma - \mu} \sigma
	- (1+\a)\intQt \m\s
	- \intQt \s\ph
	\\ && \quad \non
	+ \intQt u \sigma
	- \intQt (\pi(\ph_1) - \pi(\ph_2))\ph,
\Esist
where all the terms on the \lhs~are nonnegative. As the \rhs~is concerned,
we denote $I_1,...,I_7$ the seven integrals, in that order.
Using \eqref{young} and \eqref{propAdue} we have
\Bsist
	\non
	|I_1| = \biggl| \iot (  u + \mu + \sigma, \alpha \mu + \phi + \sigma )_* \biggr| 
	\leq \delta \int _0^t \normaVp{u + \mu + \sigma}^2
	+ \cd \int _0^t \normaVp{\alpha \mu + \phi + \sigma}^2,
\Esist
where the first term of the this inequality can be estimated by virtue of 
the embedding of $\Vp$ in $H$ and the Young~inequality as follows
\Bsist
	\non
	\d \int _0^t \normaVp{u + \mu + \sigma}^2 
	\leq c \d \iot \normaH{u + \mu\ + \sigma}^2
	\leq \frac \alpha 4  \intQt |\mu|^2 
	+ c \intQt |u|^2 
	+ c \intQt |\sigma|^2
\Esist
provided $\d$ is sufficiently small.
\dafare{Moreover, combining the \Holder~inequality
and the Sobolev continuous embedding $V \subset L^q(\Omega)$, which holds for every 
$q\in[1,6]$, we realize that
\Bsist
	\non &&
	|I_2| \leq 
	c \intQt |\phi| (|\sigma_1| + |\mu_1|) |\sigma|
	\leq c \iot \norma{\phi}_2 (\norma{\sigma_1}_4 + \norma{\mu_1}_4) \norma{\sigma}_4
	\\ \non && \quad
	\leq 
	c \iot \normaH{\phi}(\norma{\sigma_1}_V + \norma{\mu_1}_V)\normaV{\sigma}
	\leq \frac 12 \intQt \tonde{|\sigma|^2+|\nabla\sigma|^2}
	\\ \non 
	&& \quad 
	+ c \iot (\norma{\sigma_1}_V^2 + \norma{\mu_1}_V^2) \normaH \phi ^2
	\leq 
	\frac 12 \intQt \tonde{|\sigma|^2+|\nabla\sigma|^2}
	+ c \iot \normaH \phi ^2,
\Esist
where in the first line we use the \Lip\ continuity of $P$ stated by \ref{H5}, 
in the second we apply the Young inequality, while in the latter}
we made use of estimate \eqref{regstima} for the solutions $\s_1$ and $\m_1$. 
Furthermore, in view of \eqref{young}, we obtain
\Bsist
	\non
	|I_3| \leq c\intQt |\s - \m|\, |\s |
	\leq \frac \a 4 \intQt |\m|^2 
	+ c \intQt |\s|^2,
\Esist
and finally, from \eqref{young}, \eqref{stimasep}, and \ref{H7} we have that
\Bsist
	\non
	|I_4|+|I_5|+|I_6|+|I_7| \leq 
	\frac \a 4 \I2 \m
	+ \frac 12 \intQt |u|^2
	+ c \intQt |\s|^2
	+ c \intQt |\ph|^2.
\Esist
Combining the above estimates, we have shown that for every $t\in [0,T]$ it 
holds that
\Bsist
	\non
	&& \frac 12 \normaVp{(\a\m +\ph +\s)(t)}^2
	+ \frac \a 4 \intQt |\m|^2
	+ \I2 \s
	+ \frac \b 2 \iO |\ph(t)|^2 
	\\ && \quad \non
	+ \intQt |\nabla\ph|^2
	+ \frac 12 \iO |\s(t)|^2
	+ \frac 12 \intQt |\nabla \s|^2
	\leq 
	c \I2 \s
	\\ && \quad \non
	+ c \I2 \ph
	+ c \iot \normaVp{(\a\m +\ph + \s)(s)}^2 \,ds
	+ c \I2 u.
\Esist
Therefore, we invoke the Gronwall lemma and achieve
\Bsist
	&& \non
	\norma{\a\m + \ph + \s}_{\L\infty\Vp}
	+ \norma\m_{\L2H}
	+ \norma\ph_{\L\infty H \cap \L2 V}
	\\ && \quad \non
	+ \norma \s _{\L\infty H \cap \L2 V}
	\leq c \norma u _{\L2 H},
\Esist
where the variables are defined by \eqref{notdiff}. \qed

We conclude this section by proving a sharper estimate.

\proof[Proof of Theorem \ref{THcontdep2}] 
Here, we account for \eqref{contdepcontrol} and make heavily 
use of the second part of Theorem \ref{THMstate}.
We consider again that the variables are defined by \eqref{notdiff}.

{\bf First estimate:}
We consider again the system \EQContDep.
We add to both sides of \eqref{FVdifseconda} the term $-\ph$,
test \eqref{FVdiprimaprova} by $\m$, and this new second equation by $-\dt \ph$. 
Adding the equations and integrating over $Q_t$, we obtain
\Bsist
	&& \non
	\frac \a2 \IO2 \m
	+ \IN2 \m
	+ \b \I2 {\dt\ph}
	+ \frac 12 \IO2 \ph
	+ \frac 12 \INO2 \ph
	\\ && \quad \non
	\leq
	\intQt (R_1-R_2)\m
	- \intQt (F'(\ph_1)-F'(\ph_2)) \ \dt\ph
	+ \intQt \ph \,\dt\ph.
\Esist
As before we call $I_1,I_2,I_3$ the three contributions on the 
\rhs~and proceed with a separate investigation. Due to \eqref{young} 
and \Holder's~inequality, we have that
\Bsist
	\non
	&& |I_2|+|I_3| \leq 
	2\d \I2 {\dt\ph} 
	+ \cd \I2 \ph  
	+\cd \I2 {F'(\ph_1)-F'(\ph_2)}
	\\ && \quad \non
	\leq 
	2\d \I2 {\dt\ph} 
	+ \cd \I2 \ph,
\Esist
where in the last estimate we invoke the fact that, by $(ii)$ of Theorem \ref{THMstate},
$F'$ turns out to be \Lip~continuous.
Furthermore, by virtue of \ref{H5}, \eqref{contdepcontrol}, and 
\Holder's inequality and Sobolev's embeddings, we conclude
\Bsist
	\non &&
	|I_1| \leq 	
	\intQt P(\ph_2)(\s-\m)\m
	+ \intQt (P(\ph_1)-P(\ph_2))(\s_1-\m_1)\m
	\\ &&\quad \non
	\dafare{\leq 
	c \I2 \s 
	+ c \I2 \m 
	+ c \intQt |\ph|(|\s_1| + |\m_1|)\,|\m|}
	\\ &&\quad \non
	\dafare{
	\leq
	c \I2 \s 
	+ c \I2 \m 
	+c \iot \norma{\ph}_4(\norma{\s_1}_4 + \norma{\m_1}_4)\norma{\m}_2}
	\\ \non && \quad
	\dafare{
	\leq
	c \I2 \s 
	+ c \I2 \m 
	+ c \iot \normaV{\ph}^2(\normaV{\s_1}^2 + \normaV{\m_1}^2)
	\leq c \norma {u}_{\L2 H}^2,}
\Esist
\dafare{
where the fact that $\s_1$ and $\m_1$ satisfy \eqref{regstima} turn out to be
fundamental.}
On account of the previous estimates, we can choose 
$0< \d < \b/2$, and apply
the Gronwall lemma in order to conclude that
\Bsist
	\non &&
	\norma{\m}_{\L\infty H \cap \L2V}
	+ \norma{\ph}_{\H1 H \cap \L\infty V}
	\leq c \norma u _{\L2 H}.\qed
\Esist

\section{The control problem}
\label{SEC_NECESSARY}
\setcounter{equation}{0}
The current section represents the most challenging part of
the work, since it contains the proof of the \Frechet~differentiability of
the control-to-state mapping $\S$, the investigation of
both the linearized and the adjoint systems, and the necessary conditions that
a control has to satisfy to be optimal.

\subsection{Existence of optimal control}
\label{EXISTENCECONTROL}
In the following, we are going to prove the existence of an 
optimal control. We remind that in general nothing can be said about
the uniqueness. The strategy of the proof is 
quite standard and mainly lies on the semicontinuity property
of the cost functional $\cal J$ and on standard weak compactness arguments.

\proof[Proof of Theorem \ref{THexistenceofcontrol}] \quad \\
Let $\graffe {u_n}_n$ be a minimizing sequence for the control problem
\CP~constituted of elements
of $\Uad$ and for every $n \in \mathbb{N}$, let $(\m_n,\ph_n,\s_n)$ be the 
corresponding state. Therefore, the estimate \eqref{regstima} yields
that there exist $\bar{u} \in \Uad$ and a triple $(\bm,\bph,\bs)$ such that,
possibly for a subsequence which is not relabelled, it holds true the following
\Bsist
	&& \non
	u_n \to \bar{u} \ \ \hbox{weakly star in } L^\infty(Q),
	\\ \non &&
	\m_n \to \bm  \ \ \hbox{weakly star in }  \H1 H \cap \L\infty V \cap \L2 W \cap L^{\infty}(Q),
	\\ \non && 
	\ph_n \to \bph \ \ \hbox{weakly star in } W^{1,\infty}(0,T; H) \cap \H1 V \cap \L\infty W,
	\\ \non &&
	\s_n \to \bs \ \ \hbox{weakly star in }  \H1 H \cap \L\infty V \cap \L2 W.
\Esist
Furthermore, owing to standard compactness results 
\cit{(cf., e.g., \cite[Sec. 8, Cor. 4]{Simon})}, we recover even some strong convergences.
Indeed, we infer that
\Beq
	\non
	\ph_n \to \bph \ \ \hbox{strongly in } \C0 {C^0(\bar{\Omega})}.
\Eeq
This latter, paired with \ref{H5}-\ref{H8} and Theorem \ref{THMstate}, allows us
to manage the nonlinearities, since now
\Beq
	\non
	F'(\ph_n)\to F'(\bph) \
	\hbox{ and } \ P(\ph_n) \to P(\bph),
\Eeq
with the same uniform convergence.
Then, we can pass to the limit as $n$ goes to infinity in 
the variational formulation of \EQ~written for $(\m_n,\ph_n,\s_n)$. 
Therefore, we realize that $\S (\bar{u}) = (\bm,\bph,\bs)$ and
$\bar{u}$ itself are admissible solution for the \CP.
By the weak sequentially lower semicontinuity of $\cal J$ we finally realize that 
$\bar{u}$ is an optimal control that we were looking for. \qed

\subsection{Towards necessary conditions: the linearized problem}
\label{LINEARIZEDPROBLEM}
Our first efforts are intended to establish the well-posedness
of the linearized system \EQLin, namely to prove Theorem \ref{THMlinear}.
\proof[Proof of Theorem \ref{THMlinear}] \quad \\
{\it Existence:}
The \wk~spectral property of the operator $\cal A$ allow us to
apply a Faedo-Galerkin scheme. We consider the family $\graffe {w_j}_j$ of 
eigenfunctions for the eigenvalue problem
\Beq
	\non
	-\Delta w_j + w_j = \lam_j w_j \ \hbox{ in } \Omega ,
	\quad \dn w_j = 0 \ \hbox{ on } \Gamma,
\Eeq
which constitutes a Galerkin basis in $V$. Moreover, let
$\graffe {w_j}_j$ represent a complete orthonormal system in 
$(H,(\cdot,\cdot))$ which is also orthogonal in $(V,(\cdot,\cdot))$.
For fixed $n$, we set ${\cal W}_n:= \span \graffe{w_1,...,w_n}$, and 
we expect that the solutions to the approximated problem possess the following structure
\Bsist
	\non
	\et_n (x,t) = \sum_{k=1}^{n} a_k^n(t)w_k(x),
	\quad
	\th_n (x,t) = \sum_{k=1}^{n} b_k^n(t)w_k(x),
	\quad
	\r_n (x,t) = \sum_{k=1}^{n} c_k^n(t)w_k(x),
\Esist
for suitable unknown sequences $a_k^n, b_k^n, c_k^n$. Namely, we try to solve 
\EQLin~in which the variables are replaced by the above expressions and 
we will refer to this problem
as $(P_n)$. Since \eqref{EQLinseconda} only depends on
the variables $a_i^n$ and $b_i^n$, $1 \leq i \leq n$, by comparison, we can
express the unknowns $a_i^n$ in terms of $\graffe{b_1^n,...,b_n^n}$.
In this way, $(P_n)$ can be reformulated 
as a Cauchy problem for a linear system of $2n$ first-order ODE 
in the $2n$ unknowns $b_i^n, c_i^n$, $1 \leq i \leq n$.
By Cauchy-\Lip~theorem, there exists a unique solution to
this linear system satisfying
$(b_1^n,...,b_n^n, c_1^n,...,c_n^n)\in (C^{1}(0,T))^{2n}$.
This proves the existence and uniqueness of solution to $(P_n)$, and
it is \sfw~to realize that $(\et_n,\th_n, \r_n) \in {C^{1}([0,T];{\calW_n}})^3$.

At this point, we would like to obtain an existence result, 
for the solution to \EQLin~itself. To do that, we look for some a priori estimates on 
the approximated solutions that involve constants that may depend on the data 
of the problem, but are independent of $n$, thus we will be able to pass to the limit as 
$n\nearrow+\infty$ to prove the existence of solutions.
To prevent a heavy notation in the following estimates, we 
avoid writing every time the subscript $n$ under the variables, while we will
reintroduce the correct notation at the end of each estimate.

{\bf First estimate:} 
First of all, we add the term to both the members of \eqref{EQLinseconda}~$\th$. Then
we test \eqref{EQLinprima} by $\et$, this new second equation by $- \dt\th$, and
\eqref{EQLinterza} by $\r$, add the resulting equalities and integrate over $Q_t$
and by parts to obtain
\Bsist
 	\non &&
	\frac \a2 \iO |\et(t)|^2
	+ \I2  {\nabla\et}
	+ \b \I2 {\dt\th}
	+ \frac 12 \iO |\th(t)|^2
	+ \frac 12 \iO |\nabla\th(t)|^2
	\\ && \quad \non
	+ \frac 12 \iO |\r(t)|^2
	+ \I2 {\nabla\r}
	+ \intQt P(\bph)(\r - \et)^2
	=
	\intQt h\r
	\\ && \quad \non
	- \intQt F''(\bph)\,\th\, \dt\th
	+ \intQt P'(\bph)(\bs - \bm)\th(\et - \r)
	+ \intQt \th\dt\th
	\\ && \quad \non
	\leq |I_1| + |I_2| +|I_3| +|I_4|,
\Esist
where $I_1,...,I_4$ represent, in that order, the integrals in the \rhs.
It is worth to note that all the terms of the \lhs~are nonnegative since
they all contain squares and $P$ attains nonnegative values by \ref{H5}.
Clearly, by Young's inequality it turns out that
\Bsist
	\non
	|I_1| \leq \frac 12 \intQt (|h|^2+|\r|^2) ,
	\aand
	|I_2| + |I_4| \leq 2\d \I2 {\dt\th } + \cd \I2 \th,
\Esist
respectively. Moreover, by virtue of \eqref{regstima}, \eqref{stimasep},
\Holder's inequality, and the Sobolev embeddings, we have that
\dafare{
\Bsist
	\non
	&& |I_3| 
	\leq
	\intQt (|\bs|+|\bm|) \, |\th| \, (|\et|+ |\r|)
	\leq c \iot \tonde{\norma{\bs}_6 + \norma{\bm}_6} 
	\norma{\th}_3 \tonde{\norma{\et}_2 + \norma{\r}_2}
	\\ \non && \quad 
	\leq c \iot \tonde{\norma{\bs}_V + \norma{\bm}_V} 
	\norma{\th}_V \tonde{\norma{\et}_H + \norma{\r}_H}
	\leq 
	c \iot \tonde{\norma{\bs}_V^2 + \norma{\bm}_V^2} \norma{\th}_V^2
	\\ \non && \quad
	+ c \I2 \et
	+ c \I2 \r
	\leq 
	c \iot \normaV{\th}^2
	+ c \I2 \et
	+ c \I2 \r,
\Esist
where in the second line we also apply \eqref{young}.}
Thus, by fixing $ 0< \d<\b/2$, the Gronwall lemma yields
\Bsist
	&& \non
	\norma{\et_n}_{\L\infty H \cap \L2 V}
	+ \norma{\th_n}_{\H1 H \cap \L\infty V}
	\\ && \quad
	+ \norma{\r_n}_{\L\infty H \cap \L2 V}
	\leq c \norma{h}_{\L2 H}.
	\label{stimalinprima}
\Esist

{\bf Second estimate:}
We test \eqref{EQLinseconda} by $\Delta\th$, which is perfectly admissible. 
Indeed, in our approximating scheme 
$\Delta \th$ actually stands for $\Delta \th_n$, which belongs to ${\cal W}_n$. 
Using \eqref{stimasep},
the previous estimate and the Young inequality, we deduce 
$\norma{\Delta\th_n}_{\L2 H} \leq c$. Therefore, by elliptic regularity
we realize that
\Beq
	\label{stimalinseconda}
	\norma{\th_n}_{\L2 W} \leq c.
\Eeq

{\bf Third estimate:}
We multiply \eqref{EQLinprima} by $\dt \et$, \eqref{EQLinterza} by $\dt \r$, integrate 
over $Q_t$ and by parts to get
\Bsist
	\non &&
	\a \I2 {\dt \et}
	+ \frac 12 \IO2 {\nabla \et}
	+ \I2 {\dt \r}
	+ \frac 12 \IO2 {\nabla\r}
	\\ \non && \quad
	= - \intQt \dt \th \, \dt \et
	+ \intQt P'(\bph)(\bs-\bm)\th\, \dt\et
	+ \intQt P(\bph)(\r-\et)\, \dt\et	
	\\ \non && \quad
	- \intQt P'(\bph)(\bs-\bm)\th\, \dt\r
	- \intQt P(\bph)(\r-\et)\, \dt\r
	+ \intQt h \,\dt\r,
\Esist
where we denote the six terms of the \rhs~by $I_1,...,I_6$, in that order. 
\dafare{
As $I_2$ and $I_4$ are concerned, we invoke \eqref{regstima} and \eqref{stimasep} to obtain
\Bsist
	&&  \non
	|I_2| + |I_4| 
	\leq
	c \intQt (|\bs|+|\bm|) |\th| |\dt\et|
	+ c \intQt (|\bs|+|\bm|) |\th| |\dt\r|
	\\ \non && \quad
	\leq 
	c \iot (\norma{\bs}_6+\norma{\bm}_6) \norma{\th}_3 \norma{\dt\et}_2
	+ c \iot (\norma{\bs}_6+\norma{\bm}_6) \norma{\th}_3 \norma{\dt\r}_2	
	\\ \non && \quad
	\leq 
	\d \intQt (|\dt\et|^2+|\dt\r|^2)
	+ \cd \iot (\normaV{\bs}^2 +\normaV{\bm}^2) \normaV{\th}^2,
\Esist
where in the second line we apply first the \Holder~inequality and then the 
Sobolev continuous embedding of $ V \subset L^6(\Omega) $. Furthermore, the last estimate
is obtained by the Young inequality combining the fact that $\bs$ and $\bm$, as solutions, 
satisfy \eqref{regstima} and the above estimate \eqref{stimalinprima}.}
Accounting for the Young inequality, \eqref{regstima} and \eqref{stimasep},
we also conclude that 
\Bsist
	\non &&
	|I_1|
	+|I_3| 
	+ |I_5|+ |I_6|
	\leq
	2\d \intQt (|\dt\et|^2 + |\dt\r|^2)
	+\cd \intQt |\dt \th|^2	
	\\ \non && \quad
	+\cd \intQt |P(\bph)(\r-\et)|^2	
	+ \cd \intQt |h|^2	
	\\ \non && \quad
	\leq
	2\d \intQt (|\dt\et|^2 + |\dt\r|^2)
	+\cd \intQt (|\dt\th|^2 
	+ |\r|^2 + |\et|^2 + |h|^2).
\Esist
Choosing $0< \d< \min\, \graffe{\a/3,1/3}$, we can apply the Gronwall 
lemma which gives
\Beq
	\label{stimalinterza}
	\norma{\et_n}_{\H1 H \cap \L\infty V}
	+ \norma{\r_n}_{\H1 H \cap \L\infty V}
	\leq c.
\Eeq

{\bf Fourth estimate:}
Now, we test \eqref{EQLinprima} by  $-\Delta \et $, and
\eqref{EQLinterza} by $-\Delta \r$, respectively. 
Summing the resulting equalities and integrating over $Q_t$, we obtain
\Bsist
	\non &&
	\I2 {\Delta \et}
	+ \I2 {\Delta \r}
	=
	\a \intQt \dt\et \, \Delta\et
	+ \intQt \dt\th\,{\Delta \et}
	- \intQt P'(\bph)(\bs-\bm)\th\, {\Delta \et}
	\\ && \quad \non
	- \intQt P(\bph)(\r-\et) {\Delta \et}
	+ \intQt \dt\r \, {\Delta \r}
	+ \intQt P'(\bph)(\bs-\bm)\th {\Delta \r}
	\\ && \quad \non
	+ \intQt P(\bph)(\r-\et) {\Delta \r}
	- \intQt h \,{\Delta \r},
\Esist
\dafare{
where we convey to denote the integrals on the \rhs~by $I_1,...,I_8$.
Except $|I_3|$ and $|I_6|$, the other terms}
on the \rhs~that multiply $-\Delta \et $ or
$-\Delta \r$ can be easily managed by means of the Young inequality
since they are estimated with respect to the $\L2 H$ norm.
\dafare{
Moreover, owing to the \Holder~and Young inequalities and 
\eqref{stimasep}, we infer that
\Bsist
	\non && 
	|I_3|+|I_6| 
	\leq
	\intQt \Big| P'(\bph)(\bs-\bm)\th\, {\Delta \et} \, \Big|
	+ \intQt \Big| P'(\bph)(\bs-\bm)\th {\Delta \r} \, \Big|
	\\ \non && \quad
	\leq
	c \intQt (|\bs|+|\bm|) |\th| |\Delta \et|
	+ c \intQt (|\bs|+|\bm|) |\th| |\Delta \r|
	\\ \non && \quad
	\leq 
	c \iot (\norma{\bs}_6+\norma{\bm}_6) \norma{\th}_3 \norma{\Delta \et}_2
	+ c \iot (\norma{\bs}_6+\norma{\bm}_6) \norma{\th}_3 \norma{\Delta \r}_2	
	\\ \non && \quad
	\leq 
	\d \intQt (|\Delta \et|^2+|\Delta \r|^2)
	+ \cd \iot (\normaV{\bs}^2 +\normaV{\bm}^2) \normaV{\th}^2,
\Esist
where in the last two lines, we have used the Sobolev embeddings, the fact that 
$\bs$ and $\bm$ solve \EQ~and the previous estimate.
In conclusion, owing to \eqref{young} we can manage the other terms and obtain
}
\Bsist
	\non &&
	\I2 {\Delta \et}
	+ \I2 {\Delta \r}
	\leq 
	4\d \intQt ({|\Delta \et|^2+|\Delta \r|^2})
	+ \cd.
\Esist
Furthermore, we fix $0<\d<1/4$ in order to find that
 \Beq
	\label{stimalinquarta}
	\norma{\et_n}_{\L2 W}
	+ \norma{\r_n}_{\L2 W}
	\leq c.
\Eeq

{\bf Conclusion of the proof:}
Collecting all these informations, by standard compactness arguments, 
it follows that, up to a subsequence, 
suitably relabeled, $(\et_n,\th_n,\r_n)$ converges weakly star to a limit $(\et,\th,\r)$ that
solves \EQLin~and has the following regularity
\Bsist
	\non
	\et,\th,\r \in \H1 H \cap \L\infty V \cap \L2 W.
\Esist
Finally, the standard embedding results applied to each variable,
imply that they all belong to $\C0 {\Lx{r}}$ for every $r<6$.

{\it Uniqueness:} As the uniqueness is concerned we consider
\EQLin~written for the variables $(\et_i,\th_i,\r_i)$, $i=1,2$, and 
subtract the equations. Then we denote 
$\et:= \et_1-\et_2, \ \th:=\th_1-\th_2, \ \r:=\r_1-\r_2$
and observe that they solve \EQLin~with $h \equiv 0$.
Then it immediately follows that $\et=\th=\r=0$. \qed

\subsection{\Frechet~differentiability of the control-to-state mapping}
\label{FRECHET}
In the following, we prove Theorem \ref{THMFrechet}.
Let us fix $\bar{u} \in \UR$ and denote $(\bm,\bph,\bs)= \S (\bar{u})$ 
the corresponding solution to \EQ. 
Since we are going to work with small increments $h$ and $\UR$ is open, we assume $h$ to be
small enough in order that $\bar{u} + h$ belongs to $\UR$ as well. 
For $h$ fixed, we define
\Bsist
  & (\m^h,\ph^h,\s^h) := \S(\bar{u} + h ),
  \non
  \\
  & \z := \m^h - \bm - \et , \quad
  \ps := \ph^h - \bph - \th, 
  \aand
  \ch := \s^h - \bs - \r \,.
  \non
\Esist
Therefore, we aim at providing a property such as
\Beq
   \non
   \calS (\bar{u}+h) = \calS (\bar{u}) + [D{\S}(\bar{u})](h) + o(\norma{h}_{\L2 H})
   \quad{} \hbox{as} \quad{} \norma h_{\L2 H} \to 0.
\Eeq
In view of the investigation of the linearized system, by
rearranging the terms, we realize that it suffices to prove that 
\Beq
  \label{tesifrechet}
  \norma{(\z,\ps,\ch)}_{\calY}
  \leq c \norma h_{\L2 H}^2 \,
  \quad{} \hbox{as} \quad{} \norma h_{\L2 H} \to 0,
\Eeq
where $\cal Y$ stands for the space to which belongs $(\z,\ps,\ch)$.
According to Theorem \ref{THMstate} and Theorem \ref{THMlinear}, we have that
\Bsist
	\non
	{\calY}= \biggl( \H1 H \cap \L\infty V \cap \L2 W\biggr) ^3 .
\Esist

\proof[Proof of Theorem \ref{THMFrechet}]\quad \\ 
Consider \EQ~associated to $\bar{u}+h $, and subtract \EQ~associated to $\bar{u}$ and
\EQLin. By combining them, we obtain that $(\z,\ps,\ch)$ solves the following system
\Bsist
  & \a \dt \z + \dt \ps - \Delta \z = \Theta
  \quad \hbox{in $\, Q$}
  \label{EQFreprima}
  \\
  & \z = \b \dt \ps - \Delta \ps + Z
  \label{EQFreseconda}
  \quad \hbox{in $\,Q$}
  \\
  & \dt \ch - \Delta \ch = - \Theta
  \label{EQFreterza}
  \quad \hbox{in $\,Q$}
  \\[0.1cm]
  & \dn \z = \dn \ps =\dn \ch = 0
  \quad \hbox{on $\,\Sigma$}
  \label{BCEQFr}
  \\
  & \z(0)= \ps(0)= \ch(0)=0
  \quad \hbox{in $\,\Omega.$}
  \label{ICEQFr}
\Esist
\Accorpa\EQFre EQFreprima EQFreterza
where $Z$ and $\Theta$ are defined as follows
\Bsist
	\non
	& Z := F'(\ph^h)- F'(\bph)- F''(\bph)\,\th,
	\\ & \non
	\Theta := P(\ph^h)(\s^h - \m^h) - P(\bph)(\bs-\bm) - P'(\bph)(\bs -\bm)\,\th - P(\bph)(\r -\et).
\Esist
Taylor's theorem with integral remainder, and some easy calculations, allow us to 
write
\Bsist
	\non
	&&
	Z = F''(\bph)\ps + R^h_{1} (\ph^h-\bph)^2,
	\\ && \non
	\Theta = P(\bph)(\ch - \z)
	+ \tonde{P(\ph^h)- P(\bph)}\tonde{(\s^h-\bs) - (\m^h-\bm)}
	\\ && \quad \non
	+ P'(\bph)(\bs-\bm)\ps
	+ (\bs-\bm)R^h_{2}(\ph^h-\bph)^2,
\Esist
where
\Bsist
	\non
	R^h_{1}:= \int_0^1 (1-z)\,F'''(\bph + z \,(\ph^h-\bph)) dz,
	\quad 
	\quad R^h_{2}:= \int_0^1 (1-z)\,P''(\bph + z \,(\ph^h-\bph)) dz,
\Esist
respectively.
Before starting with the core of the proof, we introduce some preparatory
estimates that will be useful later on.

{\bf Preliminary estimates:}
First of all, thanks to \eqref{stimasep} and \ref{H5}-\ref{H8}, we have
\Bsist
	\label{StimaFrepropuno}
	\norma{R^h_1}_{L^{\infty} (Q)}
	+ \norma{R^h_2}_{L^{\infty} (Q)}
	\leq c.
\Esist
By \eqref{contdepcontroldue}, the previous estimate and the Sobolev embeddings, we infer 
that for every $t\in[0,T]$, it holds
\Bsist
	\non &&
	\label{StimaFrepropdue}
	\iot \Big\| {R^h_1(s)\tonde{\ph^h(s) - \bph(s)}^2}\Big\|_H^2 \,ds
	\leq c \intQt |\ph^h-\bph|^4 
	\\ && \quad
	\leq
	c \iot \norma{\ph^h-\bph}_{4}^4
	\leq c \norma{\ph^h-\bph}^4_{\L\infty V}
	\leq c \norma{h}_{\L2 H}^4.
\Esist
Furthermore, owing to \eqref{regstima}, \eqref{contdepcontrol}, \eqref{contdepcontroldue},
\Holder's~inequality, and \ref{H5}, we get
\Bsist
	\label{StimaFreproptre}
	\non &&	
	\iot \Big\| \tonde{P(\ph^h)- P(\bph)}\tonde{(\s^h-\bs) - (\m^h-\m)} \Big\|_H^2 	
	\\ \non && \quad
	\leq 
	\dafare{
	c \intQt |\ph^h-\bph|^2 (|\s^h-\bs|^2 + |\m^h-\m|^2)	
	}
	\\ && \quad \non 
	\dafare{
	\leq 
	c \iot \norma{\ph^h(s)-\bph(s)}^2_4 (\norma{\s^h(s)-\bs(s)}^2_4 + \norma{\m^h(s)-\m(s)}^2_4) \, ds	
	}
	\\ \non && \quad
	\leq 
	c \iot \normaV{\ph^h-\bph}^2 \tonde{\normaV{\s^h-\bs}^2 + \normaV{\m^h-\m}^2} 
	\leq c \,\norma{h}_{\L2 H}^4,
\Esist
\dafare{where in the third line we have applied the Sobolev embedding of $V \subset L^4(\Omega)$.}
Moreover,  from \ref{H5}, \eqref{regstima} and \eqref{stimasep}, we obtain
\Bsist
	\label{StimaFrepropcinque}
	\non
	&& \iot \Big\|{P'(\bph)(\bs-\bm)\ps}\Big\|_H^2 
	\leq 
	c \intQt \tonde{|\bs|^2+|\bm|^2}|\ps|^2
	\\ && \quad
	\leq
	c \iot \tonde{\normaV{\bs}^2 + \normaV{\bm}^2}\normaV{\ps}^2
	\leq c \iot \normaV{\ps}^2.
\Esist
Finally, thanks to \eqref{StimaFrepropuno}, \Holder's~inequality, \eqref{regstima}, 
\eqref{contdepcontrol}, \eqref{contdepcontroldue}, and to the Sobolev embeddings,
we have
\dafare{
\Bsist
	\label{StimaFrepropquattro}
	\non &&
	\iot \Big\|{(\bs-\bm)R^h_{2}(\ph^h-\bph)^2}\Big\|_H^2
	\leq c \intQt (|\bs|^2 + |\bm|^2) |{\ph^h-\bph}|^{4}
	\\ && \quad \non
	\leq c \iot \tonde{\norma{\bs(s)}^2_6 + \norma{\bm(s)}^2_6} \norma{\ph^h(s)-\bph(s)}^{4}_6 \,ds
	\\ && \quad
	\leq c \iot \tonde{\normaV{\bs}^2 + \normaV{\bm}^2} \normaV{\ph^h-\bph}^{4}
	\leq c \norma{h}_{\L2 H}^4.
\Esist
}
Now, we start with the actual estimates.

{\bf First estimate:} 
First, we add to both sides of \eqref{EQFreseconda} the term $ \ps$, then
we multiply \eqref{EQFreprima} by $\z$, this new second equation by $-\dt\ps$,
and \eqref{EQFreterza} by $\ch$. Adding the resulting equations and integrating over $Q_t$, we get
\Bsist
	\non
	&& \frac \a2 \iO |\z(t)|^2
	+ \I2 {\nabla\z}
	+ \frac 12 \iO |\ps(t)|^2
	+ \frac 12 \iO |\nabla \ps(t)|^2
	+ \b \I2 {\dt\ps}
	\\ && \quad \non
	+ \iO |\ch(t)|^2
	+ \I2 {\nabla\ch}
	= 
	\intQt \Theta \z
	- \intQt F''(\bph)\,\ps\,\dt\ps
	\\ && \quad \non
	- \intQt R^h_1(\ph^h-\bph)^2\,\dt\ps
	+ \intQt \ps \, \dt\ps
	- \intQt \Theta \ch,
\Esist
where the last five integrals of the \rhs~are denoted by
$I_1,...,I_5$, in this order. 
Simply using \eqref{stimasep}, \eqref{young}, and \eqref{StimaFrepropdue}, we deduce
\Bsist
	\non
	&& |I_2|+|I_3|+|I_4| \leq 
	3 \d \I2 {\dt\ps}
	+ \cd \I2 \ps
	+ \cd \I2 {R^h_1(\ph^h-\bph)^2} 
	\\ \non
	&& \quad \leq 
	3\d \I2 {\dt\ps}
	+ \cd \I2 \ps
	+ \cd \norma{h}_{\L2 H}^4.
\Esist
Moreover, we have
\Bsist
	\non
	&& |I_1| \leq \
	\bigr| 
		\intQt
		P(\bph)(\ch - \z)\, \z 
		+ \tonde{P(\ph^h)- P(\bph)}\tonde{(\s^h-\bs) - (\m^h-\m)}\z		
		\\ && \quad \non
		+ P'(\bph)(\bs-\bm)\, \ps \, \z
		+ (\bs-\bm)R^h_{2}(\ph^h-\bph)^2 \, \z
	\bigl|
	\\ && \quad \non
	\leq
	c \I2 {\ch}
	+ c \I2 {\z}
	+ c \iot \normaV{\ps}^2
	+ c \norma{h}_{\L2 H}^4 ,
\Esist
owing to the Young inequality, \ref{H5}, \eqref{regstima}, 
\eqref{stimasep}, and the estimates \accorpa{StimaFreproptre}{StimaFrepropquattro}. 
The last term $I_5$ is treated the same way, while it is referred to the variable $\ch$
instead of $\z$. 
Choosing $0< \d <\b/3$, we can apply the Gronwall lemma in order to realize that
\Bsist
	&& \non
	\norma{\z}_{\L\infty H \cap \L2 V}
	+ \norma{\ps}_{\H1 H \cap \L\infty V }
	\\ && \quad
	+ \, \norma{\ch}_{\L\infty H \cap \L2 V}
	\leq
	c \norma{h}_{\L2 H}^2.
	\label{Frechetprimastima}
\Esist

{\bf Second estimate:} 
Accounting for the previous estimate,
by comparison in \eqref{EQFreseconda}, we easily conclude that
\Beq
	\non
	\norma{\Delta \ps}_{\L2 H} \leq c \norma{h}_{\L2 H}^2.
\Eeq

\dafare{
{\bf Third estimate:} 
To recover the stated regularity, let us reformulate the equations 
\eqref{EQFreprima} and \eqref{EQFreterza} as follows
\Beq
	\non
	\a \dt \z -\Delta \z =  \Theta - \dt\ph := g_1, 
	\aand 
	\dt\ch -\Delta \ch = \Theta := g_2.
\Eeq
Accounting for \eqref{Frechetprimastima}, we realize that both the forcing
terms $g_1$ and $g_2$ have been already estimated in $\L2 H$.
Moreover, owing to the smoothness of the initial conditions \eqref{BCEQFr},
the parabolic regularity theory (see, e.g., \cit{\cite{LDY}}) gives
\Beq
	\non
	\non
	\norma{\z}_{\H1 H \cap \L\infty V \cap \L2 W} 
	+ \norma{\chi}_{\H1 H \cap \L\infty V \cap \L2 W} 
	\leq c \norma{h}_{\L2 H}^2.
\Eeq
This proves \eqref{tesifrechet}, that is the \Frechet~differentiability of $\S$. \qed}

%

At this point Corollary \ref{CORprimanec} immediately follows from \eqref{abstrnec} by
direct calculations.

\subsection{Adjoint problem}
\setcounter{equation}{0}
The last part of our work regards the improvement of \eqref{primanec} 
by dealing with the system \EQAgg. In fact, our aim is to prove Theorem \ref{THMadjoint}.
\proof[Proof of Theorem \ref{THMadjoint}]\quad \\ 
{\it Existence:}
As in the proof of Theorem \ref{THMlinear}, we apply a Faedo-Galerkin
scheme based on a basis $\graffe{{w_j}}_{j}\subset W$, and we again refer to ${\cal W}_n$ as to
the space generated by the first $n$ eigenvectors.
We look for approximated solutions of the form
\Bsist
	\non
	q_n (x,t) = \sum_{k=1}^{n} a_k^n(t)w_k(x),
	\quad
	p_n (x,t) = \sum_{k=1}^{n} b_k^n(t)w_k(x),
	\quad
	r_n (x,t) = \sum_{k=1}^{n} c_k^n(t)w_k(x),
\Esist
which satisfies, $\aat$ the following problem
\Bsist	
	\non
	\label{Faedouno}
	&& \b(\dt q_n, v)	
	+ (-\dt p_n, v)
	{-} (\nabla q_n, \nabla v)
	- (F''(\bph)q_n,v)
	\\ && \quad
	+ (P'(\bph)(\bs-\bm)(r_n-p_n), v)
	= (\bQ (\bph-\phQ), v)
	\quad \hbox{for all $v \in {\cal W}_n$,}
	\\
	\label{Faedodue}
	&& (q_n,v)
	{-\a}(\dt p_n, v)
	{+} (\nabla p_n, \nabla v)
	+(P(\bph)(p_n-r_n), v) = 0
	\quad \hbox{for all $v \in {\cal W}_n$,}
	\\
	\non
	&& (-\dt r_n,v)
	+ (\nabla r_n, \nabla v)
	+ (P(\bph)(r_n-p_n), v)
	\\ && \quad
	\label{Faedotre}
	= \andrea{(\bQh (\bs-\sQ), v)}
	\quad \hbox{for all $v \in {\cal W}_n$,}
	\\ &&  \non
	p_n(T) -\b q_{n}(T) = \mathbb{P}\bigr(\bO(\bph(T)-\phO)\bigl), \quad
	\a  p_n(T) = 0,
	\\ && \quad
	\quad r_n (T) = \andrea{\mathbb{P}\bigr(\bOh(\bs(T)-\sO)\bigl)},	
	\label{Faedoquattro}
\Esist
\Accorpa\EQFaedo Faedouno Faedoquattro
where $\mathbb{P}$ represents the orthogonal projection in $H$ onto ${\cal W}_n$.
Arguing as before, we can easily conclude that the backward-in-time problem
\EQFaedo~admits a unique solution triple that satisfies the following regularity
$(q_n, p_n, r_n) \in {\bigl(\W{1,\infty}{\calW_n}}\bigr)^3$.
So, to ensure the existence of the adjoint problem, we need to provide some a 
priori estimates independent of $n$ in order to apply standard compactness arguments 
and motivate rigorously the passage to the limit as $n\nearrow+\infty$.

As for the notation, we again adopt the convention used in the proof of Theorem \ref{THMlinear}.

{\bf First estimate:} 
First, we add to both sides of \eqref{Faedodue} the term $p$. 
Then, we test \eqref{Faedouno} by $-q$, this new second equation
by $- \dt p $, and \eqref{Faedotre} by $r$. Finally, we add these equations and 
integrate over $\Omega \times [t, T] = : Q_t^T$ and by parts to find the following identity
\Bsist
	&& \non
	\frac \b 2 \iO|q(t)|^2
	+ \intQtT \dt p \, q
	+ \intQtT |{\nabla q}|^2
	- \intQtT \dt p \, q
	+ \a \intQtT |{\dt p}|^2
	\\ && \quad \non
	+ \frac 12 \iO |\nabla p(t)|^2
	+ \frac 12 \iO |p(t)|^2
	+ \frac 12 \iO |r(t)|^2
	+ \intQtT |{\nabla r}|^2
	\\ && \quad \non
	=
	\andrea{
	\frac 12 \iO|r(T)|^2
	+ \frac \b 2 \iO|q(T)|^2}
	+ \frac 12 \iO |\nabla p(T)|^2
	+ \frac 12 \iO |p(T)|^2
	\\ && \quad \non
	+ \intQtT P'(\bph)(\bs-\bm)(r-p)\, q
	- \intQtT F''(\bph)q^2
	- \intQtT \bQ(\bph-\phQ)q
	\\ && \quad \non
	\andrea{+ \intQtT \bQh (\bs-\sQ)r}
	- \intQtT P(\bph)(r-p)r
	- \intQtT P(\bph)(r-p)\, \dt p
	- \intQtT p \, \dt p.
\Esist
Let us note that two terms cancel out and that the first four 
integrals of the \rhs~can be explicitly written using \eqref{Faedoquattro} and are bounded due to
\ref{H1}-\ref{H2}. Let us call,
in the order, $I_1,...,I_7$ the other terms. Using \eqref{stimasep} and \ref{H2}, we have
\Beq
	\non
	|I_2|+|I_3| + |I_4| \leq c + c \intQtT |q|^2 + c \intQtT |r|^2.
\Eeq
In addition, \eqref{stimasep} and \eqref{young} yield that 
\Bsist
	\non &&
	|I_5|+|I_6| + |I_7| \leq 
	2\d \intQtT |{\dt p}|^2
	+c \intQtT |r|^2
	+ \cd \intQtT |{P(\bph)(r-p)}|^2
	\\ \non && \quad
	+ \cd \intQtT | p|^2
	\leq 
	2\d \intQtT | {\dt p}|^2
	+ \cd \intQtT | r |^2
	+ \cd \intQtT |p|^2.
\Esist
Finally, we obtain from \ref{H5}, \eqref{regstima}, \eqref{stimasep},
\eqref{young}, the Sobolev embeddings, and the \Holder~inequality that
\dafare{
\Bsist
	\non
	&&
	|I_1| \leq
	c \IT2 {P'(\bph)(\bs-\bm)q}
	+ c \IT2 {r-p}	
	\\ \non && \quad
	\leq 
	c \intQtT(|\bs|^2+|\bm^2|)|q||q|
	+ c \IT2 r 
	+ c \IT2 p
	\\ \non && \quad
	\leq 
	c \int_t^T (\norma{\bs}_6^2 + \norma{\bm}_6^2)\norma{q}_6 \norma{q}_2
	+ c \IT2 r 
	+ c \IT2 p
	\\ \non && \quad
	\leq 
	\frac 12 \intQtT (|q|^ 2+ |\nabla q|^2)
	+ c \int_t^T (\normaV{\bs}^4 + \normaV{\bm}^4)\norma{q}_2^2
	+ c \IT2 r 	
	+ c \IT2 p
	\\ \non && \quad
	\leq
	\frac 12 \IT2 {\nabla q}
	+c \IT2 r + c \IT2 p+ c \IT2 q.
\Esist
}
We now fix $0< \d <\a/2$, and applying the backward in time Gronwall lemma, we infer that
\Bsist
	\non &&
	\norma {q_n}_{\L\infty H \cap \L2 V}
	+ \norma {p_n}_{\H1 H \cap \L\infty V}
	+ \norma {r_n}_{\L\infty H \cap \L2 V }
	\leq c.
\Esist

{\bf Second estimate:}
We test \eqref{Faedodue} by $\Delta p$.
Using the Young inequality and the previous estimate it is quite easy to realize that 
\Beq
	\non
	\norma{\Delta p_n }_{\L2 H} \leq c,
\Eeq
whence, \dafare{from elliptic regularity,} we infer that
\Beq
	\non
	\norma{p_n}_{\L2 W} \leq c.
\Eeq

{\bf Third estimate:} 
We now test \eqref{Faedouno} by $\dt q$. Integrating over $\Omega \times [t,T]$ and
by parts, we obtain that  
\Bsist
	\non &&
	\b \IT2 {\dt q}
	+ \frac 12 \IO2 {\nabla q}
	=
	\frac 12 \iO |\nabla q(T)|^2
	+ \intQtT \dt p \, \dt q
	\\ \non && \quad
	- \intQtT P'(\bph)(\bs-\bm)(r-p)\,\dt q
	+ \intQtT F''(\bph)q\,\dt q
	+ \intQtT \bQ(\bph- \phQ) \, \dt q,
\Esist
and we denote by $I_1,...,I_4$ the last four summands on the \rhs.
Note that the first term on the \rhs~is finite by \eqref{Faedoquattro} and \ref{H2}.
A simple application of \eqref{stimasep} and of the Young inequality show that
\Beq
	\non
	|I_1| 
	+ |I_3| 
	+ |I_4|
	\leq  
	\cd
	+ 3 \d \IT2 {\dt q} 
	+ \cd \IT2 {\dt p}
	+ \cd \IT2 {q}.
\Eeq
\dafare{
Furthermore, owing to the Young inequality and to the Sobolev embeddings, we also have that
\Bsist
 	&&  \non
 	|I_2|
 	\leq 
	c \intQtT (|\bs|+|\bm|)(|r|+|p|)|\dt q|
	\leq 
	c \int_t^T (\norma{\bs}_6 + \norma{\bm}_6)(\norma{r}_3+\norma{p}_3 )\norma{\dt q}_2
	\\ \non && \quad
	\leq
	\d \IT2 {\dt q}
	+ \cd \int_t^T  (\norma{\bs}_V^2 + \norma{\bm}_V^2)(\norma{r}_V^2+\norma{p}_V^2 ),
\Esist
where all the terms on the \rhs~of both these inequalities have been already estimated above.
Therefore, fixing $0< \d <\b/4$, we conclude}
\Beq
	\non
	\norma{q_n}_{\H1 H \cap \L\infty V} 
	\leq c.
\Eeq

{\bf Fourth estimate:}
Arguing exactly as above, by testing \eqref{Faedotre} by $-\dt r$, we also infer that 
\Beq
	\non
	\norma{r_n}_{\H1 H \cap \L\infty V} \leq c.
\Eeq


{\bf Fifth estimate:}
Moreover, accounting for the above estimates and the Young inequality,
by taking $-\Delta r$ and $\Delta q$ as test functions in \eqref{EQAggprima} 
and \eqref{EQAggterza}, respectively, we can easily deduce that 
\Beq
	\non
	\norma{q_n}_{\L2 W} + \norma{r_n}_{\L2 W} \leq c.
\Eeq

{\bf Conclusion of the proof:}
It follows from the above a priori estimates that there exist
functions $(q,p,r)$ such that, possibly for some subsequence which is again indexed by $n$,
the following convergences
\Bsist
	\non
	q_n \to q , \quad p_n \to p, \quad r_n \to r \quad \hbox{weakly star in } \H1 H \cap \L\infty V \cap \L2 W
\Esist
hold.
Moreover, by continuous embedding, also
\Bsist
	\non
	q_n \to q , \quad p_n \to p, \quad  r_n \to r  \quad \hbox{weakly in } \C0 V.
\Esist

It is now a standard matter to verify that $(q,p,r)$ is in fact a solution 
to the system \EQAgg~satisfying \eqref{regadj}.
\\
{\it Uniqueness:} As before we denote
$q:=q_1-q_2, \ p:=p_1-p_2, \ r:= r_1-r_2$, where $(q_i,p_i,r_i)$, $i=1,2$, are two solutions
to \EQAgg. If we consider the system obtained by subtracting the corresponding 
equations each others, we can repeat
the argument of the existence and realize that $q=p=r=0$.
\qed

\subsection{Final necessary condition}
We are now in the position to eliminate $\th$ and $\r$ from \eqref{primanec}.
This procedure automatically leads to \eqref{secondanec} and prove Theorem \ref{THMsecondanec}.
\dafare{
\proof[Proof to Theorem \ref{THMsecondanec}]
Comparing \eqref{primanec} with \eqref{secondanec}, we realize that it sufficies to show that 
\Bsist
	\non
	&&  \intQ  r h 
	=  \bQ \intQ (\bph - \phQ)\th 
	+ \bO \iO (\bph(T) - \phO)\th(T)
	\\ && \quad 
	+ \andrea{\bQh \intQ (\bs - \sQ)\r }
	+ \andrea{\bOh \iO (\bs(T) - \sO)\r(T)},
\Esist
where $\th$ and $\r$ solve the linearized system \EQLin~with $h =v- \bar u$.
Indeed, if this equality are satisfied \eqref{secondanec} directly follows by \eqref{primanec} by
a mere substitution.
In this direction, owing to \EQLin, we have that the following equalities are satisfied:
\Bsist
	\non
	&& 0= \intQ q \, [\et - \b \dt \th + \Delta \th - F''(\bph)\,\th] ,
	\\  && \non
	0 = \intQ p \,[\a \, \dt \et + \dt\th - \Delta \et 
	- P'(\bph)(\bs-\bm)\, \th -P(\bph)(\r- \et)],
	\\ && \non
	0 = \intQ r \, [\dt \r - \Delta \r + P'(\bph)(\bs-\bm)\, \th + P(\bph)(\r-\et) - h] .
\Esist
Hence, summing the above equalities and integrating by parts, we realize that
\Bsist
	&& \non
	0 = \intQ q \et + \b \intQ  \dt q\, \th - \b \iO \th(T)q(T) 
	+ \intQ  \Delta q \,\th - \intQ F''(\bph)\, q  \, \th
	- \a \intQ \dt p \,\et 
	\\ && \non \quad 
	+ \a \iO p(T)\et(T) - \intQ \dt p \, \th + \iO p(T)\th(T) - \intQ \Delta p \, \et
	- \intQ P'(\bph)(\bs-\bm)\, p \, \th 
	\\ && \non \quad
	- \intQ P(\bph) p \, \r + \intQ P(\bph) p \, \et
	- \intQ \dt r \, \r + \iO r(T)\r(T) - \intQ \Delta r \, \r 
	\\ && \non \quad  
	+ \intQ P'(\bph)(\bs-\bm) r\, \th 
	+ \intQ P(\bph) r \, \r - \intQ P(\bph) r \, \et - \intQ r \, h,
\Esist
where, after the time integration, only the final conditions are remained
since the initial value of the linearized variables are all zero by \eqref{ICEQLin}.
Moreover, in the integration by parts of the terms with the Laplacian, we also account for the
homogeneous Neumann boundary conditions \eqref{BCEQLin}. 
Therefore, by rearranging the above equality we get
\Bsist
	&& \non
	\intQ r \, h = \intQ  [q -\a \dt p - \Delta p + P(\bph)(p -r)] \, \et
	\\ && \non \quad  
	+ \intQ [\b \dt q - \dt p + \Delta q - F''(\bph)q + P'(\bph)(\bs - \bm)(r-p)] \, \th 
	\\ && \non \quad  
	+ \intQ [-\dt r - \Delta r + P(\bph)(r - p)] \, \r
	\\ && \non \quad  
	+\iO [- \b \th(T)q(T) + \a \et(T)p(T) + p(T) \th(T)  + r(T)\r(T)],
\Esist
and invoking the adjoint system \EQAgg, this latter reduces to
\Bsist
	\non
	&&  \intQ  r h 
	=  \bQ \intQ (\bph - \phQ)\th 
	+ \bO \iO (\bph(T) - \phO)\th(T)
	\\ && \quad 
	+ \andrea{\bQh \intQ (\bs - \sQ)\r }
	+ \andrea{\bOh \iO (\bs(T) - \sO)\r(T)},
\Esist
which is the equality we were looking for.
\qed
}

\Brem
\label{REM_modcost} 
Let us slightly digress to point out a mathematical issue. 
The choice of the tracking type cost functional \eqref{costfunct} is 
essentially led by the model interpretation. Indeed, from a 
mathematical point of view, only little rearrangements are needed to treat the
more general version
\Bsist
	\non
	\hat{\J} (\ph,\m,\s, u)  : = 
	\J (\ph,\s,u)	
	+\frac {b_5} 2 \norma{\m - \mQ}_{L^2(Q)}^2
	+\frac {b_6} 2 \norma{\m(T)-\mO}_{L^2(\Omega)}^2,
\Esist
in which all the variables appear.
At this stage, we understood the natural requirements
on the constants and on the targets that are necessary to give sense to these lines.
As the necessary condition is concerned, we expect something like \eqref{primanec} in 
which the following additional terms on its \lhs
\Bsist
	\non
   b_5 \intQ (\bm - \mQ)\et  
  + b_6\iO (\bm(T) - \mO)\et(T)
\Esist
occur.
Moreover, the adjoint system will read exactly as \EQAgg, but instead of 
\eqref{EQAggseconda} and \eqref{BCEQAgg} we should have
\Bsist
  \non
  &q -\a \dt p - \Delta p + P(\bph)(p -r) = {b_5}(\bm - \mQ)
  \quad \hbox{in $\,Q$,}
  \\ \non
  & \quad \hbox{ and } \quad 
  \a p(T) = {b_6}(\bm(T) - \mO) 
  \quad \hbox{in $\Omega$,}
\Esist
respectively. About the existence result, note that
the presence of this new term on the \rhs~of \eqref{EQAggseconda} does not add difficulties 
since it can be easily handled by the Young inequality.
In fact, only \sfw~modifications are needed to extend the proof of Theorem \ref{THMadjoint}
to this general framework. \dafare{In a similar way also the new final 
condition can be handled.}
\Erem

To conclude, let us mention that for
forthcoming contributions, it will be interesting to couple our
study for the control problem \CP~with asymptotic analysis as $\a$ and $\b$ go to zero. Of course, this
would require less generality for the potentials, in order to handle the 
passage to the limit, as pointed out in \cite{CGH}, \cite{CGRS_ASY} and \cite{CGRS_VAN}.

\subsection*{Acknowledgments}
The author wishes to thanks Professor Pierluigi Colli for several useful 
comments and suggestions without which this paper cannot be possible. 
Moreover, the author would like to thank especially one of the referees 
for the careful reading and the precious suggestions which 
have improved the manuscript.

\vspace{3truemm}

\Begin{thebibliography}{10}

\footnotesize

\bibitem{BRZ}
H. Brezis,
``Op\'erateurs maximaux monotones et semi-groupes de contractions dans les
espaces de Hilbert'', North-Holland Math. Stud. {\bf 5}, North-Holland, Amsterdam, 1973.

\bibitem{CGH}
P. Colli, G. Gilardi and D. Hilhorst,
On a Cahn-Hilliard type phase field system related to tumor growth,
{\it Discrete Contin. Dyn. Syst.} {\bf 35} (2015), 2423-2442.

\bibitem{CGMR_cons}
P. Colli, G. Gilardi, G. Marinoschi and E. Rocca,
Optimal control for a conserved phase field system with a possibly singular potential,
{\it Evol. Equ. Control Theory\/} {\bf 7} (2018), 95-116.

\bibitem{CGMR_sing}
P. Colli, G. Gilardi, G. Marinoschi and E. Rocca,
Optimal control for a phase field system with a possibly singular potential,
{\it Math. Control Relat. Fields\/} {\bf 6} (2016), 95-112.

\bibitem{CGRS_OPT}
P. Colli, G. Gilardi, E. Rocca and J. Sprekels,
Optimal distributed control of a diffuse interface model of tumor growth,
{\it Nonlinearity} {\bf 30} (2017), 2518-2546.

\bibitem{CGRS_ASY}
P. Colli, G. Gilardi, E. Rocca and J. Sprekels,
Asymptotic analyses and error estimates for a \CH~type phase field system modeling tumor growth,
{\it Discrete Contin. Dyn. Syst. Ser. S} {\bf 10} (2017), 37-54.

\bibitem{CGRS_VAN}
P. Colli, G. Gilardi, E. Rocca and J. Sprekels,
Vanishing viscosities and error estimate for a Cahn–Hilliard type phase field system related to tumor growth,
{\it Nonlinear Anal. Real World Appl.} {\bf 26} (2015), 93-108.

\bibitem{CGS_nonst}
P. Colli, G. Gilardi and J. Sprekels,
Optimal boundary control of a nonstandard viscous \CH~system with dynamic boundary condition,
{\it Nonlinear Anal.\/} {\bf 170} (2018), 171-196.

\bibitem{CGS14}
P. Colli, G. Gilardi and J. Sprekels,
Optimal velocity control of a viscous Cahn--Hilliard system with convection and dynamic boundary conditions, 
{\it SIAM J. Control Optim.\/} {\bf 56} (2018), 1665-1691.

\bibitem{CGS_OPT}
P. Colli, G. Gilardi and J. Sprekels,
A boundary control problem for the viscous \CH~equation 
with dynamic boundary conditions,
{\it Appl. Math. Opt.\/} {\bf 73} (2016) 195-225.

\bibitem{CGS}
P. Colli, G. Gilardi and J. Sprekels,
On the \CH~equation with dynamic 
boundary conditions and a dominating boundary potential,
{\it J. Math. Anal. Appl.\/} {\bf 419} (2014) 972-994.

\bibitem{CS}
P. Colli and J. Sprekels,
Optimal control of an Allen--Cahn equation 
with singular potentials and dynamic boundary condition,
{\it SIAM J. Control Optim.\/} {\bf 53} (2015) 213-234.

\bibitem{DFRGM}
M. Dai, E. Feireisl, E. Rocca, G. Schimperna, M. Schonbek,
Analysis of a diffuse interface model of multispecies tumor growth,
{\it Nonlinearity\/} {\bf  30} (2017), 1639.

\bibitem{FGR}
S. Frigeri, M. Grasselli, E. Rocca,
On a diffuse interface model of tumor growth,
{\it  European J. Appl. Math.\/} {\bf 26 } (2015), 215-243. 

\bibitem{FLRS}
S. Frigeri, K.F. Lam, E. Rocca, G. Schimperna,
On a multi-species Cahn--Hilliard-Darcy tumor growth model with singular potentials,
preprint arXiv:1709.01469 (2017), 1-41, {\it Comm Math Sci.\/}, to appear (2018). 

\dafare{
\bibitem{GARL_1}
H. Garcke and K. F. Lam,
Well-posedness of a Cahn-–Hilliard–-Darcy system modelling tumour
growth with chemotaxis and active transport,
{\it European. J. Appl. Math.} {\bf 28 (2)} (2017), 284–316.
\bibitem{GARL_2}
H. Garcke and K. F. Lam,
Analysis of a Cahn--Hilliard system with non--zero Dirichlet 
conditions modeling tumor growth with chemotaxis,
{\it Discrete Contin. Dyn. Syst.} {\bf 37 (8)} (2017), 4277-4308.
\bibitem{GARL_3}
H. Garcke and K. F. Lam,
Global weak solutions and asymptotic limits of a Cahn--Hilliard--Darcy system modelling tumour growth,
{\it AIMS Mathematics} {\bf 1 (3)} (2016), 318-360.
\bibitem{GARL_4}
H. Garcke and K. F. Lam,
On a Cahn--Hilliard--Darcy system for tumour growth with solution dependent source terms, 
in {\sl Trends on Applications of Mathematics to Mechanics}, 
E.~Rocca, U.~Stefanelli, L.~Truskinovski, A.~Visintin~(ed.), 
{\it Springer INdAM Series} {\bf 27}, Springer, Cham, 2018, 243-264.
\bibitem{GARLR}
H. Garcke, K. F. Lam and E. Rocca,
Optimal control of treatment time in a diffuse interface model of tumor growth,
{\it Appl. Math. Optim.} {\bf 28} (2017), {1-50}.
\bibitem{GAR}
H. Garcke, K. F. Lam, R. N\"urnberg and E. Sitka,
A multiphase Cahn--Hilliard--Darcy model for tumour growth with necrosis,
{\it Mathematical Models and Methods in Applied Sciences} {\bf 28 (3)} (2018), 525-577.
}

\bibitem{GiMiSchi} 
G. Gilardi, A. Miranville and G. Schimperna,
On the \CH~equation with irregular potentials and dynamic boundary conditions,
{\it Commun. Pure Appl. Anal.\/} {\bf 8} (2009) 881-912.

\bibitem{HDZO}
A. Hawkins-Daruud, K. G. van der Zee and J. T. Oden, Numerical simulation of
a thermodynamically consistent four-species tumor growth model, Int. J. Numer.
{\it Math. Biomed. Engng.} {\bf 28} (2011), 3–24.

\bibitem{HKNZ}
D. Hilhorst, J. Kampmann, T. N. Nguyen and K. G. van der Zee, Formal asymptotic
limit of a diffuse-interface tumor-growth model, 
{\it Math. Models Methods Appl. Sci.} {\bf 25} (2015), 1011–1043.

\bibitem{LDY}
O.A. Lady\v zenskaja, V.A. Solonnikov and N.N. Uralceva,
``Linear and quasilinear equations of parabolic type'',
Mathematical Monographs Volume {\bf 23},
{\it American mathematical society}, Providence, 1968.

\bibitem{Lions}
J.-L. Lions,
``\'Equations diff\'erentielles op\'erationnelles et probl\`emes aux limites'',
Grundlehren, Band~111,
Springer-Verlag, Berlin, 1961.

\bibitem{LioMag}
J.-L. Lions and E. Magenes,
``Non-homogeneous boundary value problems and applications'',
Vol.~I,
Springer, Berlin, 1972.

\bibitem{MirZelik}
A. Miranville and S. Zelik,
Robust exponential attractors for \CH~type equations with singular potentials,
{\it Math. Methods Appl. Sci.\/} {\bf 27} (2004) 545--582.

\bibitem{Simon}
J. Simon,
{Compact sets in the space $L^p(0,T; B)$},
{\it Ann. Mat. Pura Appl.~(4)\/} 
{\bf 146} (1987) 65-96.

\bibitem{Trol}
F. Tr\"oltzsch,
{``Optimal control of partial differential equations:
theory, methods and applications''},
American mathematical society, Providence, 2010.
Originally published in German, under the title: ``Optimale Steuerung partieller
Differentialgleichungen'', 2005.

\bibitem{WZZ}
X. Wu, G.J. van Zwieten and K.G. van der Zee, Stabilized second-order splitting
schemes for \CH~models with applications to diffuse-interface tumor-growth models, 
{\it Int. J. Numer. Meth. Biomed. Engng.} {\bf 30} (2014), 180-203.

\End{thebibliography}

\End{document}

\bye